\documentclass[11pt]{article}

\usepackage{amsfonts}
\usepackage{enumitem}
\usepackage{shuffle,yfonts,amsmath,amssymb,amsthm}
\usepackage{hyperref}
\hypersetup{pdfstartview={XYZ null null 1.00}}
\usepackage{pdfpages}

 \allowdisplaybreaks

\catcode`!=11
\let\!int\int \def\int{\displaystyle\!int}
\let\!lim\lim \def\lim{\displaystyle\!lim}
\let\!sum\sum \def\sum{\displaystyle\!sum}
\let\!sup\sup \def\sup{\displaystyle\!sup}
\let\!inf\inf \def\inf{\displaystyle\!inf}
\let\!cap\cap \def\cap{\displaystyle\!cap}
\let\!max\max \def\max{\displaystyle\!max}
\let\!min\min \def\min{\displaystyle\!min}
\let\!frac\frac \def\frac{\displaystyle\!frac}
\catcode`!=12

\let\oldsection\section
\renewcommand\section{\setcounter{equation}{0}\oldsection}

\allowdisplaybreaks

\def\R{\mathbb{R}}
\def\N{\mathbb{N}}
\def\Z{\mathbb{Z}}
\def\Q{\mathbb{Q}}
\def\CC{\mathbb{C}}
\def\sn{{\rm sn}}
\def\cn{{\rm cn}}
\def\dn{{\rm dn}}
\def\su{\sum\limits_{n=1}^\infty}

\newtheorem{defn}{Definition}[section]
\newtheorem{thm}{Theorem}[section]
\newtheorem{lem}[thm]{Lemma}
\newtheorem{cor}[thm]{Corollary}
\newtheorem{con}[thm]{Conjecture}
\newtheorem{re}{Remark}[section]
\newtheorem{exa}{Example}[section]
\DeclareMathOperator{\Res}{Res}

\begin{document}
\title {\bf Reciprocal Hyperbolic Series of Ramanujan Type}
\author{
{Ce Xu${}^{a,}$\thanks{Email: cexu2020@ahnu.edu.cn, first author, ORCID 0000-0002-0059-7420.}\ and Jianqiang Zhao${}^{b,}$\thanks{Email: zhaoj@ihes.fr, corresponding author,  ORCID 0000-0003-1407-4230.}}\\[1mm]
\small a. School of Mathematics and Statistics, Anhui Normal University, Wuhu 241002, PRC\\
\small b. Department of Mathematics, The Bishop's School, La Jolla, CA 92037, USA}

\date{}
\maketitle \noindent{\bf Abstract} This paper presents an approach to summing a few families of infinite series involving hyperbolic functions, some of which were first studied by Ramanujan. The key idea is based on their contour integral representations and residue computations with the help of some well-known results of Eisenstein series given by Ramanujan, Berndt et al. As our main results, several series involving hyperbolic functions are evaluated and expressed in terms of $z={}_2F_1(1/2,1/2;1;x)$ and $z'=dz/dx$. When a certain parameter in these series is equal to $\pi$ the series are expressed in closed forms in terms of some special values of the Gamma function. Moreover, many new illustrative examples are presented.



\section{Introduction}
In the Entry 16(x) of Chapter 17 in his second notebook \cite{B1991,R2012}, Ramanujan wrote down the identity
\begin{equation}\label{equ:Rama1}
\sum\limits_{n=1}^\infty \frac {(-1)^{n-1}(2n-1)^2}{\cosh((2n-1)\pi/2)} = \frac{{{\pi ^{3/2}}}}
{{2\sqrt 2 {\Gamma ^6}\left( {\frac{3}
{4}} \right)}}.
\end{equation}
He also evaluated a few similar series some of which involved hyperbolic sine function $\sinh$.
Our main goal of this paper is to extend and study systematically more of such series as defined in Definition \ref{de1}.

Let $\N$ be the set of positive integers, $\N_0:=\N\cup \{0\}$, $\Z$ the ring of integers, $\Q$ the field of rational numbers, $\R$ the field of real numbers, and $\CC$ the field of complex numbers.

For any $a\in\CC$ and $n\in\N_0$, let $(a)_n$ denote the ascending Pochhammer symbol defined by the following:
\begin{align*}
(a)_0:=1\quad {\rm and}\quad (a)_n:=\frac {\Gamma(a+n)}{\Gamma(a)}=a(a+1)\cdots(a+n-1),
\end{align*}
where $\Gamma(z)$ denotes the Gamma function. When $\Re(z) > 0$
\[\Gamma(z) := \int\limits_0^\infty  {{e^{ - t}}{t^{z - 1}}dt}.\]

The Gaussian or ordinary hypergeometric function ${_2}F_1(a,b;c;x)$ is defined for $|x|<1$ by the power series
\begin{align}
{_2}F_1(a,b;c;x)=\sum\limits_{n=0}^\infty \frac{(a)_n(b)_n}{(c)_n} \frac {x^n}{n!}\quad (a,b,c\in\mathbb{C}).
\end{align}
Next, we recall that the complete elliptic integral of the first kind is defined by (cf. \cite{WW1966})
\begin{align}
K=K(x)=K(k^2):=\int\limits_{0}^{\pi/2}\frac {d\varphi}{\sqrt{1-k^2\sin^2\varphi}}=\frac {\pi}{2} {_2}F_{1}\left(\frac {1}{2},\frac {1}{2};1;k^2\right).\label{1.3}
\end{align}
Here $x=k^2$ and $k\ (0<k<1)$ is the modulus of $K$. The complementary modulus $k'$ is defined by $k'=\sqrt {1-k^2}$. Furthermore,
\begin{align*}
K':=K(k'^2)=\int\limits_{0}^{\pi/2}\frac {d\varphi}{\sqrt{1-k'^2\sin^2\varphi}}=\frac {\pi}{2} {_2}F_{1}\left(\frac {1}{2},\frac {1}{2};1;1-k^2\right).
\end{align*}
Similarly, the complete elliptic integral of the second kind is denoted by (cf. \cite{WW1966})
\begin{align}
E=E(x)=E(k^2):=\int\limits_{0}^{\pi/2} \sqrt{1-k^2\sin^2\varphi}d\varphi=\frac {\pi}{2} {_2}F_{1}\left(-\frac {1}{2},\frac {1}{2};1;k^2\right).\label{1.5}
\end{align}

In order to better state our main results, we shall henceforth adopt the notations of Ramanujan (see Berndt's book \cite{B1991}). Let
\begin{align}\label{1.6}
x:=k^2,\quad y=y(x):=\pi \frac {K'}{K}, \quad q=q(x):=e^{-y},\quad z=z(x):=\frac {2}{\pi}K,\quad z':=\frac{dz}{dx}.
\end{align}
Using the identity $(a)_{n+1}=a(a+1)_n$, one sees that (see \cite[p. 405, 16.3.1]{AskeyOD1999})
\begin{align}
\frac {d^n}{dx^n}{_2}F_1(a,b;c;x)=\frac{(a)_n(b)_n}{(c)_n}{_2}F_1(a+n,b+n;c+n;x)\quad(n\in\N_0).\label{1.7}
\end{align}
Then,
\begin{align*}
\frac {d^nz}{dx^n}=\frac{\left(1/2\right)^2_n}{n!}{_2}F_1\left(\frac{1}{2}+n,\frac{1}{2}+n;1+n;x\right).
\end{align*}
Applying the identity (see Theorem 3.5.4(i) in Section 3.5 of Chapter 3 of Andrews, Askey and Roy book \cite{A2000})
\begin{align*}
{_2}F_1\left(a,b;\frac {a+b+1}{2};\frac 1{2}\right)=\frac{\Gamma\left(\frac{1}{2}\right)\Gamma\left(\frac{a+b+1}{2}\right)}
{\Gamma\left(\frac{a+1}{2}\right)\Gamma\left(\frac{b+1}{2}\right)},
\end{align*}
we get
\begin{align*}
\left.\frac {d^nz}{dx^n}\right|_{x=1/2}=\frac{(1/2)^2_n\sqrt{\pi}}{\Gamma^2\left(\frac{n}{2}+\frac {3}{4}\right)}.
\end{align*}
Setting $n=0,1,2$ in the above equation yields
$$
z\Big(\frac 1{2}\Big)=\frac{\Gamma^2(1/4)}{2\pi^{3/2}},\ z'\Big(\frac 1{2}\Big)= \frac {4\sqrt{\pi}}{\Gamma^2(1/4)}\quad {\rm and}\quad z''\Big(\frac 1{2}\Big)=\frac{\Gamma^2(1/4)}{2\pi^{3/2}},
$$
where we have used the following two relations of the Gamma function:
$$
\Gamma(x+1)=x\Gamma(x)\quad{\rm and}\quad \Gamma(x)\Gamma(1-x)=\frac{\pi}{\sin(\pi x)}.
$$

We now define a few Ramanujan type sums involving reciprocal hyperbolic functions.
\begin{defn}\label{de1} Let $m\in\N$ and $p\in\Z$. Define
\begin{alignat*}{4}
&S_{p,m}(y):=\sum\limits_{n=1}^\infty \frac {n^p}{\sinh^m(ny)},\quad &&{\bar S}_{p,m}(y):=\sum\limits_{n=1}^\infty \frac {n^p}{\sinh^m(ny)}(-1)^{n-1},\\
&C_{p,m}(y):=\sum\limits_{n=1}^\infty \frac {n^p}{\cosh^m(ny)},\quad &&{\bar C}_{p,m}(y):=\sum\limits_{n=1}^\infty \frac {n^p}{\cosh^m(ny)}(-1)^{n-1},\\
&S'_{p,m}(y):=\sum\limits_{n=1}^\infty \frac {(2n-1)^p}{\sinh^m((2n-1)y/2)},\quad &&{\bar S}'_{p,m}(y):=\sum\limits_{n=1}^\infty \frac {(2n-1)^p}{\sinh^m((2n-1)y/2)}(-1)^{n-1},\\
&C'_{p,m}(y):=\sum\limits_{n=1}^\infty \frac {(2n-1)^p}{\cosh^m((2n-1)y/2)},\quad &&{\bar C}'_{p,m}(y):=\sum\limits_{n=1}^\infty \frac {(2n-1)^p}{\cosh^m((2n-1)y/2)}(-1)^{n-1}.
\end{alignat*}
\end{defn}

\begin{re}
For the convenience of the reader, here is the mnemonic scheme we have used:
\begin{center}
$S \leftrightarrow \sinh$, $C \leftrightarrow \cosh$, an apostrophe $' \leftrightarrow$ odd numbers, a bar $\overline{\phantom{B}} \leftrightarrow$ alternating series.
\end{center}
\end{re}

By \eqref{1.6}, the eight series in Definition \ref{de1} can be rewritten as Eisenstein series
\begin{alignat*}{4}
&S_{p,m}(y)=2^m \su \frac{n^pq^{mn}}{(1-q^{2n})^m},\quad && {\bar S}_{p,m}(y)=2^m \su \frac{n^pq^{mn}}{(1-q^{2n})^m}(-1)^{n-1},\\
&C_{p,m}(y)=2^m\su \frac{n^p q^{mn}}{(1+q^{2n})^m},\quad &&{\bar C}_{p,m}(y)=2^m\su \frac{n^p q^{mn}}{(1+q^{2n})^m}(-1)^{n-1},\\
&S'_{p,m}(y)=2^m \su \frac {(2n-1)^p q^{(n-1/2)m}}{(1-q^{2n-1})^m},\quad &&{\bar S}'_{p,m}(y)=2^m \su \frac {(2n-1)^p q^{(n-1/2)m}}{(1-q^{2n-1})^m}(-1)^{n-1},\\
&C'_{p,m}(y)=2^m \su \frac {(2n-1)^p q^{(n-1/2)m}}{(1+q^{2n-1})^m},\quad &&{\bar C}'_{p,m}(y)=2^m \su \frac {(2n-1)^p q^{(n-1/2)m}}{(1+q^{2n-1})^m}(-1)^{n-1}.
\end{alignat*}

Infinite series involving hyperbolic functions have attracted the attention of many authors. In particular, Ramanujan evaluated many such series in his notebooks \cite{R2012} and lost notebook \cite{R1988}.
For example, Ramanujan obtained \eqref{equ:Rama1} in his second notebook \cite{B1991,R2012}.
Further, he showed in the Entries 15, 16 and 17 of Chapter 17 in Berndt's book \cite{B1991} that one can also evaluate the following more general sums in closed forms
\begin{align*}
S_{2p-1,1}(y),\ C_{2p-2,1}(y),\ S'_{2p-1,1}(y),\quad C'_{2p-2,1}(y),\ {\bar C}'_{2p-1,1}(y),
\end{align*}
for all $p\in\N$. These evaluations are in terms of $z={}_2F_1(1/2,1/2;1,k^2)$. In \cite{L1974,L1975,L1978},
Ling used Weierstrassian elliptic functions to obtain a certain closed forms for the fifteen series
\begin{align*}
S_{0,2m}(y),\ {\bar S}_{0,2m}(y),\ C_{0,m}(y),\ {\bar C}_{0,m}(y),\ S'_{0,2m}(y),\ {\bar S}'_{0,2m-1}(y),\ C'_{0,m}(y),
\end{align*}
and
\begin{align*}
S_{2p-1,1}(y),\ {\bar S}_{2p-1,1}(y),\ C_{2p,1}(y),\ {\bar C}_{2p,1}(y),\ S'_{2p-1,1}(y),\ {\bar S}'_{2p,1}(y),\ C'_{2p,1}(y),\ {\bar C}'_{2p-1,1}(y),
\end{align*}
for all $m,p\in\N$, and summed the series in terms of Gamma functions for $y=\pi,\ \sqrt{3}\pi$ and ${\pi}/{\sqrt{3}}$. Further, Ling \cite{L1988} extended the summation of the above series to $y=\sqrt{n}\pi$ and $\pi/\sqrt{n}$ for all $n\in\N$.
Zucker \cite{Z1979} showed that the series $C_{0,m}(y)$ and $C'_{0,m}(y)$ can be expressed in closed forms for all $m\in\Z$. Moreover, he stated that many other series of hyperbolic functions may be expressed in closed form, e.g.
\[\sum\limits_{n=1}^\infty \frac{n^p}{\sinh((2n-1)y/2)},\ \sum\limits_{n=1}^\infty \frac{n^p}{\cosh((2n-1)y/2)},\quad {\rm etc.}\]
Zucker also studied several infinite series involving exponential and hyperbolic functions, which depend on a certain parameter $y/\pi$. When $y/\pi$ is the square root of a rational number the sums may be expressed in terms of special values of the $\Gamma$-function and other well-known transcendental numbers. Further details may be found in \cite{Z1984}.

Some other related results for infinite series involving hyperbolic functions may be seen in the works of \cite{AB2009,AB2013,B1977,B1978,B2004,B2016,BS2017,C1889,L2011,S1954,Ya-2018}. For instance, Berndt \cite{B1977,B1978} found a lot of identities about infinite series using a certain modular transformation formula that originally stems from the generalized Eisenstein series. There are also a lot of recent contributions on Eisenstein series involving hyperbolic functions (see \cite{T2015,T2008,T2010,T2012}) and the references therein.

However, the evaluations of series with $m,p\geq 2$ defined in Definition \ref{de1} has not been done so far.
The purpose of this paper is to present a method to evaluate the following eight series of quadratic hyperbolic functions:
\begin{align*}
S_{2p,2}(y),\ {S}'_{2p,2}(y),\ C_{2p,2}(y),\ {C}'_{2p,2}(y),\ {\bar S}_{2p,2}(y),\ {\bar S}'_{2p-1,2}(y),\ {\bar C}_{2p,2}(y)\ {\rm and}\ {\bar C}'_{2p-1,2}(y)\quad (p\in\N).
\end{align*}
In sections \ref{sec3.1} and \ref{sec3.2}, we use the residue theorem and asymptotic formulas of trigonometric and hyperbolic functions at the poles to establish four equations involving infinite series of trigonometric and hyperbolic functions (see Theorem \ref{thm3.1}). Then we apply these four equations and several well-known results given by Ramanujan, Berndt et al. to prove that the four alternating sums ${\bar S}_{2p,2}(y),\ {\bar S}'_{2p-1,2}(y),\ {\bar C}_{2p,2}(y)\ {\rm and}\ {\bar C}'_{2p-1,2}(y)$ can be expressed in terms of $z$ and $z'$ (defined by \eqref{1.6}). Furthermore, in section \ref{sec3.3}, we use the results about Eisenstein series \eqref{3.52}--\eqref{3.59} and \eqref{3.63}--\eqref{3.66} established by Ramanujan, and the transformations \eqref{2.23}--\eqref{2.25} to prove that the four series $S_{2p,2}(y),\ { S}'_{2p,2}(y),\ C_{2p,2}(y)\ {\rm and}\ { C}'_{2p,2}(y)$
can be evaluated in terms of $z$ and $z'$. In the last section, we give numerous examples of series involving hyperbolic functions. Moreover, we prove that the eight series ($p\geq k$)
\begin{align*}
S_{2p,2k}(y),\ {\bar S}_{2p,2k}(y),\ C_{2p,2k}(y),\ {\bar C}_{2p,2k}(y),\ S'_{2p,2k}(y),\ {\bar S'}_{2p-1,2k}(y),\ C'_{2p,2k}(y),\ {\bar C'}_{2p-1,2k}(y)
\end{align*}
can be evaluated by recurrence formulas.

\section{Some Lemmas and Transformations}

In this section we give some lemmas and transformations, which will be useful in the development of the main results of this paper.

\subsection{Three Lemmas}

\begin{lem}\emph{(Residue Theorem, \cite{FS1998})}\label{lem2.1}
Let $\xi (s)$ be a kernel function and let $r(s)$ be a function which is $O(s^{-2})$ at infinity. Then
\begin{align}\label{e1}
\sum\limits_{\alpha\in T} \underset{s=\alpha}\Res  \big(r(s)\xi(s)\big)   + \sum\limits_{\beta\in S} \underset{s=\beta}\Res \big(r(s)\xi(s)\big)  = 0,
\end{align}
where $S$ is the set of poles of $r(s)$ and $T$ is the set of poles of $\xi(s)$ that are not poles $r(s)$. Here $\underset{s=\alpha}\Res \big(r(s)\big) $ denotes the residue of $r(s)$ at $s=\alpha$. The kernel function $\xi(s)$ is meromorphic in the whole complex plane and satisfies $\xi (s)=o(s)$ over an infinite collection of circles $|s|= {\rho_k}$ with ${\rho _k} \to \infty$ (i.e., $\xi (s)/s\to 0$ as $|s|\to \infty$).
\end{lem}

\begin{lem} \emph{(cf. \cite{X2018})}\label{lem2.2} Let $n$ be an integer, then the following asymptotic formulas hold:
\begin{align}
&\pi \cot \left( {\pi s} \right)\mathop  = \limits^{s \to n} \frac{1}{{s - n}} - 2\sum\limits_{k = 1}^\infty  {\zeta \left( {2k} \right){{\left( {s - n} \right)}^{2k - 1}}} ,\label{e2}
\\
&\frac{\pi }
{{\sin \left( {\pi s} \right)}}\mathop  = \limits^{s \to n} {\left( { - 1} \right)^n}\left( {\frac{1}
{{s - n}} + 2\sum\limits_{k = 1}^\infty  {\bar \zeta \left( {2k} \right){{\left( {s - n} \right)}^{2k - 1}}} } \right),\label{e3}
\\
&\pi \coth \left( {\pi s} \right)\mathop  = \limits^{s \to ni} \frac{1}
{{s - ni}} - 2\sum\limits_{k = 1}^\infty  {{{\left( { - 1} \right)}^k}\zeta \left( {2k} \right){{\left( {s - ni} \right)}^{2k - 1}}} ,\label{e4}
\\
&\frac{\pi }
{{\sinh \left( {\pi s} \right)}}\mathop  = \limits^{s \to ni} {\left( { - 1} \right)^n}\left( {\frac{1}
{{s - ni}} + 2\sum\limits_{k = 1}^\infty  {{{\left( { - 1} \right)}^k}\bar \zeta \left( {2k} \right){{\left( {s - ni} \right)}^{2k - 1}}} } \right),\label{e5}
\\
& \pi \tan \left( {\pi s} \right)\mathop  = \limits^{s \to (n-1/2)}  - \frac{1}
{{s - \frac{{2n - 1}}
{2}}} + 2\sum\limits_{k = 1}^\infty  {\zeta \left( {2k} \right){{\left( {s - \frac{{2n - 1}}
{2}} \right)}^{2k - 1}}} ,\label{e6}
\\
&\pi \tanh \left( {\pi s} \right)\mathop  = \limits^{s \to (n-1/2)i} \frac{1}
{{s - \frac{{2n - 1}}
{2}i}} - 2\sum\limits_{k = 1}^\infty  {{{\left( { - 1} \right)}^k}\zeta \left( {2k} \right){{\left( {s- \frac{{2n - 1}}
{2}i} \right)}^{2k - 1}}} ,\label{e7}
\\
&\frac{\pi }
{{\cos \left( {\pi s} \right)}}\mathop  = \limits^{s \to n-1/2} {\left( { - 1} \right)^n}\left\{ {\frac{1}
{{s - \frac{{2n - 1}}
{2}}} + 2\sum\limits_{k = 1}^\infty  {\bar \zeta \left( {2k} \right){{\left( {s - \frac{{2n - 1}}
{2}} \right)}^{2k - 1}}} } \right\},\label{e8}\\
& \frac{\pi }
{{\cosh \left( {\pi s} \right)}}\mathop  = \limits^{s \to (n-1/2)i} {\left( { - 1} \right)^n}i\left\{ {\frac{1}
{{s - \frac{{2n - 1}}
{2}i}} + 2\sum\limits_{k = 1}^\infty  {{{\left( { - 1} \right)}^k}\bar \zeta \left( {2k} \right){{\left( {s - \frac{{2n - 1}}
{2}i} \right)}^{2k - 1}}} } \right\},\label{e9}
\end{align}
where $\zeta(s)$ and ${\bar \zeta}(s)$ denote the Riemann zeta function and alternating Riemann zeta function, respectively, namely,
\begin{align*}
\zeta(s):=\sum\limits_{n = 1}^\infty {\frac {1}{n^{s}}}\quad (\Re(s)>1)\quad {\rm and}\quad \bar \zeta (s) := \sum\limits_{n = 1}^\infty  {\frac{{{{\left( { - 1} \right)}^{n - 1}}}}{{{n^s}}}}\quad ( \Re(s) \ge 1).
\end{align*}
\end{lem}

When applying the above lemma, we often need to use Euler's famous formula (see \cite{A2000})
\begin{align*}
\zeta(2m) = \frac{(-1)^{m-1}B_{2m}}{2(2m)!} (2\pi)^{2m}
\end{align*}
for all $m\in\N$, where $B_{2m}$ is Bernoulli number.

\begin{lem}\label{lem2.3} \emph{(cf. \cite{B1991})} For any $p\in\N$, the four series
\begin{align*}
S_{2p-1,1}(y),\ S'_{2p-1,1}(y),\ C_{2p-2,1}(y)\ {\rm and}\ C'_{2p-2,1}(y)
\end{align*}
can be expressed in terms of $z$ defined by \eqref{1.6}.
\end{lem}
\begin{proof}
On the one hand, by theorems of Hermite, which may be found in Cayley's book \cite{C1961} or see Berndt's book \cite{B1991} and paper \cite{B2016}, for $u=2Kt/\pi$ and $|u|<K'$,
\begin{align}
\sn u=& u-(1+k^2)\frac{u^3}{3!}+(1+14k^2+k^4)\frac{u^5}{5!}-(1+135k^2+135k^4+k^6)\frac{u^7}{7!}\nonumber\\
&+(1+1228k^2+5478k^4+1228k^6+k^8)\frac{u^9}{9!}+\cdots,\label{2.11}\\
\cn u=& 1-\frac{u^2}{2!}+(1+4k^2)\frac{u^4}{4!}-(1+44k^2+16k^4)\frac{u^6}{6!}\nonumber\\
&+(1+408k^2+912k^4+64k^6)\frac{u^8}{8!}+\cdots,\label{2.12}\\
\dn u=&1-k^2\frac{u^2}{2!}+k^2(4+k^2)\frac{u^4}{4!}-k^2(16+44k^2+k^4)\frac{u^6}{6!}\nonumber\\
&+k^2(64+912k^2+408k^4+k^6)\frac{u^8}{8!}+\cdots,\label{2.13}
\end{align}
where $\sn u,\ \cn u$ and $\dn u$ denote the Jacobian elliptic functions. (The power series exansions of the functions $\sn u,\ \cn u,\ \dn u$ and related functions can be routinely found. They are given to quite high orders by Hancock \cite{H1958}).
On the other hand, by results of Jacobi, which may be found in Whittaker and Watson's treatise \cite{WW1966} (or see Guo, Tu and Wang \cite{G1989}), for $u=2Kt/\pi$ and $|u|<K'$,
\begin{align}
\sn u=&\frac{2\pi}{Kk}\sum\limits_{n=0}^\infty \frac{q^{(2n+1)/2}\sin(2n+1)t}{1-q^{2n+1}}\nonumber\\
=&\frac{2\pi}{Kk}\sum\limits_{j=0}^\infty \frac{(-1)^jt^{2j+1}}{(2j+1)!}\sum\limits_{n=0}^\infty \frac{(2n+1)^{2j+1}q^{(2n+1)/2}}{1-q^{2n+1}}\nonumber\\
=&\frac{\pi}{Kk}\sum\limits_{j=0}^\infty \frac{(-1)^j(\pi u/2K)^{2j+1}}{(2j+1)!}S'_{2j+1,1}(y),\label{2.14}\\
\cn u=&\frac{2\pi}{Kk}\sum\limits_{n=0}^\infty \frac{q^{n+1/2}\cos(2n+1)t}{1+q^{2n+1}}\nonumber\\
=&\frac{2\pi}{Kk}\sum\limits_{j=0}^\infty \frac{(-1)^jt^{2j}}{(2j)!}\sum\limits_{n=0}^\infty \frac{(2n+1)^{2j}q^{n+1/2}}{1+q^{2n+1}}\nonumber\\
=& \frac{\pi}{Kk}\sum\limits_{j=0}^\infty \frac{(-1)^j(\pi u/2K)^{2j}}{(2j)!}C'_{2j,1}(y),\label{2.15}\\
\dn u=&\frac{\pi}{2K}+\frac{2\pi}{K}\sum\limits_{n=1}^\infty \frac{q^n \cos(2nt)}{1+q^{2n}}\nonumber\\
=&\frac{\pi}{2K}+\frac{2\pi}{K}\sum\limits_{j=0}^\infty \frac{(-1)^j(2t)^{2j}}{(2j)!}\sum\limits_{n=1}^\infty \frac{n^{2j}q^n}{1+q^{2n}}\nonumber\\
=&\frac{\pi}{2K}+\frac{\pi}{K}\sum\limits_{j=0}^\infty \frac{(-1)^j(\pi u/K)^{2j}}{(2j)!}C_{2j,1}(y).\label{2.16}
\end{align}
Furthermore, from Greenhill's treatise \cite{G1892}, we have
\begin{align}
\sn^2u&=\Big(1-\frac{E}{K}\Big)\frac{1}{k^2}-\frac{\pi^2}{K^2k^2}\sum\limits_{n=1}^\infty \frac{n \cos(n\pi u/K)}{\sinh(ny)}\nonumber\\
&=\Big(1-\frac{E}{K}\Big)\frac{1}{k^2}-\frac{\pi^2}{K^2k^2}\sum\limits_{j=0}^\infty \frac{(-1)^j(\pi u/K)^{2j}}{(2j)!}S_{2j+1,1}(y),\label{2.17}
\end{align}
where $K$ and $E$ are represented by \eqref{1.3} and \eqref{1.5}, respectively.

Moreover, \eqref{2.14} yields
\begin{align}
\sn^2 u=u^2-(1+k^2)\frac{u^4}{3}+(2+13k^2+2k^4)\frac{u^6}{45}+\cdots.\label{2.18}
\end{align}
Then, combining \eqref{2.11}-\eqref{2.18} and applying the notations of \eqref{1.6}, we may deduce the desired results by comparing coefficients of $u^j$ in equations \eqref{2.11}--\eqref{2.18} on the right-hand side. The proof of the lemma is now complete.
\end{proof}

In the Entries 15, 16 and 17 of Chapter 17 of \cite{B1991}, Berndt presented explicit formulas for
\begin{align*}
&S_{2p-1,1}(y)\quad (p=2,3,4,5), \ S'_{2p-1,1}(y)\quad (p=1,2,3,4,5,6),\\
&C_{2p-2,1}(y)\quad (p=1,2,3,4,5),\ C'_{2p-2,1}(y)\quad (p=1,2,3,4,5).
\end{align*}
Note that from \cite{B1991}, we can also obtain
\begin{align*}
S_{1,1}(y)=\frac1{4} \big(xz^2-2x(1-x)zz'\big).
\end{align*}

\subsection{Three Transformations}

In Chapters 17 and 18 of his book \cite{B1991}, Berndt described some procedures in the theory of elliptic functions by which ``new" formulas can be produced from ``old" formulas. Here, we need to cite the following three implications:
\begin{align}
&{\rm Given}\quad \Omega(x,e^{-y},z)=0\ \Longrightarrow\ \Omega(1-x,e^{-\pi^2/y},yz/\pi)=0,\label{2.19}\\
&{\rm Given}\quad \Omega(x,e^{-y},z)=0\ \Longrightarrow\ \Omega\Big(\frac{x}{x-1},-e^{-y},z\sqrt{1-x}\Big)=0,\label{2.20}\\
&{\rm Given}\quad \Omega(x,e^{-y},z)=0\ \Longrightarrow\ \Omega\left(\Bigg(\frac{1-\sqrt{1-x}}{1+\sqrt{1-x}}\Bigg)^2,e^{-2y},\frac{1}{2}z(1+\sqrt{1-x})\right)=0,\label{2.21}
\end{align}
where $\Omega(x,e^{-y},z)$ is an arbitrary function of $x,e^{-y}$ and $z$.
In fact, if we add $z'$ as one component on the left-hand side of transformations \eqref{2.19}--\eqref{2.21}, namely, if we change $\Omega(x,e^{-y},z)$ to the form
\begin{align}\label{2.22}
\Omega(x,e^{-y},z,z')=0,
\end{align}
where $\Omega(x,e^{-y},z,z')$ is an arbitrary function of $x,e^{-y},z$ and $z'$, then after some elementary manipulation we can derive the following three new transformations.

\begin{thm} Set $r(x)=\sqrt{1-x}$. Given the formula \eqref{2.22}, we can deduce the formulas
\begin{align}
&\Omega\left(1-x,e^{-\pi^2/y},yz/\pi,\frac{1}{\pi}\bigg(\frac{1}{x(1-x)z}-yz'\bigg)\right)=0,\label{2.23}\\
&\Omega\left(\frac{x}{x-1},-e^{-y},z\sqrt{1-x},(1-x)^{3/2}\Big(\frac{z}{2}-(1-x)z'\Big)\right)=0,\label{2.24}\\
&\Omega\left(\Bigg(\frac{1-r(x)}{1+r(x)}\Bigg)^2,e^{-2y},\frac{1}{2}z(1+r(x)),\frac{(1+r(x))^3}{4(1-r(x))}\Big(-\frac{z}{2}+ \big(r(x)+1-x\big)z'\Big)\right)=0.\label{2.25}
\end{align}
\end{thm}
\begin{proof} Suppose that $x_1,\ y_1,\ z_1$ and $z'_1$ is another set of parameters such that
\[\Omega\big(x_1,e^{-y_1},z_1,z'_1\big)=0.\]
According to \eqref{2.19}--\eqref{2.21}, in order to prove \eqref{2.23}--\eqref{2.25} we only need to prove
\begin{align*}
&z'_1=\frac{1}{\pi}\bigg(\frac{1}{x(1-x)z}-yz'\bigg),\ {\rm if}\ x_1= 1-x,\\
&z'_1=(1-x)^{3/2}\Big(\frac{z}{2}-(1-x)z'\Big),\ {\rm if}\ x_1=\frac{x}{x-1},\\
&z'_1=\frac{(1+r(x))^3}{4(1-r(x))}\Big(-\frac{z}{2}+ \big(r(x)+1-x\big)\Big)z',\ {\rm if}\ x_1=\Bigg(\frac{1-r(x)}{1+r(x)}\Bigg)^2.
\end{align*}
By formula \eqref{1.7} and the definition of $z$, we obtain
\begin{align*}
z'=\frac1{4} {}_2F_1\bigg(\frac3{2},\frac 3{2};2;x\bigg).
\end{align*}
Hence, if $x_1=1-x$, then
\begin{align*}
z'_1=&\frac1{4} {}_2F_1\bigg(\frac3{2},\frac 3{2};2;1-x\bigg)=-\frac{d}{dx}{}_2F_1\bigg(\frac{1}{2},\frac1{2};1;1-x\bigg)\\
=& -\frac{d}{dx}(yz/\pi)=\frac{1}{\pi}\bigg(\frac{1}{x(1-x)z}-yz'\bigg);
\end{align*}
if $x_1=\frac{x}{x-1}$, then
\begin{align*}
z'_1=&\frac1{4} {}_2F_1\bigg(\frac3{2},\frac 3{2};2;x_1\bigg)=-(x-1)^2\frac{d}{dx}{}_2F_1\bigg(\frac{1}{2},\frac1{2};1;x_1\bigg)\\
=&-(x-1)^2\frac{d}{dx}(r(x)z)=(1-x)^{3/2}\Big(\frac{z}{2}-(1-x)z'\Big);
\end{align*}
and if $x_1=\Big(({1-r(x)})/({1+r(x)})\Big)^2$, then
\begin{align*}
z'_1=&\frac1{4} {}_2F_1\left(\frac3{2},\frac 3{2};2;x_1\right)\\
=&\frac{(1+r(x))^3r(x)}{2(1-r(x))}\frac{d}{dx}{}_2F_1\left(\frac1{2},\frac 1{2};1;x_1\right)\\
=&\frac{(1+r(x))^3r(x)}{4(1-r(x))}\frac{d}{dx}(z(1+r(x)))\\
=&\frac{(1+r(x))^3}{4(1-r(x))}\Big(-\frac{z}{2}+ \big(r(x)+1-x\big)z'\Big),
\end{align*}
where we have used the three relations (see \cite{B1989,B1991})
\begin{align*}
\frac{dy}{dx}=-\frac{1}{x(1-x)z^2},\quad {}_2F_1\bigg(\frac{1}{2},\frac1{2};1;\frac{x}{x-1}\bigg)=r(x)z
\end{align*}
and
$$
{}_2F_1\left(\frac1{2},\frac 1{2};1;\Bigg(\frac{1-r(x)}{1+r(x)}\Bigg)^2\right)=\frac{1}{2}z(1+r(x)).
$$
Thus, the proofs of \eqref{2.23}--\eqref{2.25} are completed.
\end{proof}

It should be emphasized that \cite{B1991} also contains many other types of transformations.

\section{Evaluations of some reciprocal hyperbolic series}

In this section, we will prove our main results by using the residue theorem, transformations \eqref{2.23}--\eqref{2.25} and some well-known results of Ramanujan, Berndt et al. First, we define a few auxiliary functions.

\begin{defn} Let $s\in\CC$ and $a,b,\theta\in\R$ with $|\theta|<2b \pi$, $a\neq 0$ and $b>0$. Define
\begin{align*}
&F_1(s,\theta;a,b):=\frac{\pi^3\cosh(\theta s)}{\sin(a \pi s)\sinh^2(b \pi s)},
\qquad H_1(s,\theta;a,b):=\frac{\pi^4\sinh(\theta s)}{\sin^2(a \pi s)\sinh^2(b \pi s)}, \\
&F_2(s,\theta;a,b):=\frac{\pi^3\sinh(\theta s)}{\cos(a \pi s)\cosh^2(b \pi s)},
\qquad H_2(s,\theta;a,b):=\frac{\pi^4\sinh(\theta s)}{\cos^2(a \pi s)\cosh^2(b \pi s)},\\
&F_3(s,\theta;a,b):=\frac{\pi^3\cosh(\theta s)}{\sin(a \pi s)\cosh^2(b \pi s)},
\qquad H_3(s,\theta;a,b):=\frac{\pi^4\sinh(\theta s)}{\sin^2(a \pi s)\cosh^2(b \pi s)},\\
&F_4(s,\theta;a,b):=\frac{\pi^3\sinh(\theta s)}{\cos(a\pi s)\sinh^2(b \pi s)}.
\end{align*}
\end{defn}

From Lemma \ref{lem2.2}, for all $n\in \Z$ we may deduce the asymptotic expansions of a few reciprocal quadratic trigonometric and hyperbolic functions as follows:
\begin{align}
\left(\frac{\pi}{\sin(a \pi s)}\right)^2\mathop  = \limits^{s \to n/a} &\, \frac{1}{(as-n)^2}+2\zeta(2)+6\zeta(4)(as-n)^2,\label{3.1}\\
\left(\frac{\pi}{\sinh(a \pi s)}\right)^2\mathop  = \limits^{s \to ni/a}&\, \frac{1}{(as-ni)^2}-2\zeta(2)+6\zeta(4)(as-ni)^2,\label{3.2}\\
\left(\frac{\pi}{\cos(a \pi s)}\right)^2\mathop  = \limits^{s \to (2n-1)/(2a)}&\, \frac 1{\Bigg(as-\frac{2n-1}{2}\Bigg)^2}+2\zeta(2)+6\zeta(4)\left(as-\frac{2n-1}{2}\right)^2,\label{3.3}\\
\left(\frac{\pi}{\cosh(a \pi s)}\right)^2\mathop  = \limits^{s \to (2n-1)i/(2a)}&\, -\frac 1{\Bigg(as-\frac{2n-1}{2}i\Bigg)^2}+2\zeta(2)-6\zeta(4)\left(as-\frac{2n-1}{2}i\right)^2.\label{3.4}
\end{align}
Similarly, many other asymptotic expansions can also be established. For example,
\begin{align*}
&\left(\frac{\pi}{\sin(a \pi s)}\right)^3\mathop  = \limits^{s \to n/a}(-1)^n\left(\frac{1}{(as-n)^3}+3\frac{\zeta(2)}{as-n}+\frac{51}{4}\zeta(4)(as-n)\right).
\end{align*}

\subsection{Some preliminary results}\label{sec3.1}

\begin{thm}\label{thm3.1} For $a,b,\theta\in\R$ with $|\theta|<2b \pi$, $a\neq 0$ and $b>0$,
\begin{align}
&b^2\pi^2\sum\limits_{n=1}^\infty \frac{\cosh(\theta n/a)}{\sinh^2(bn\pi/a)}(-1)^n+\theta \pi a \su \frac{\sin(n\theta/b)}{\sinh(an\pi/b)}\nonumber\\
&\quad+a^2\pi^2\su \frac{\cos(n\theta/b)\cosh(an\pi/b)}{\sinh^2(an\pi/b)} +\frac{\theta^2}{4}+\frac{\zeta(2)}{2}(a^2-2b^2)=0,\label{3.5}\\
&b^2\pi^2\su\frac{\sinh((2n-1)\theta/2a)}{\cosh^2((2n-1)\pi b/2a)}(-1)^n-\theta \pi a \su \frac{\cos((2n-1)\theta/2b)}{\cosh((2n-1)\pi a/2b)}\nonumber
\\&\quad+a^2\pi^2\su \frac{\sin((2n-1)\theta/2b)\sinh((2n-1)\pi a/2b)}{\cosh^2((2n-1)\pi a/2b)}=0,\label{3.6}\\
&b^2\pi^2\su \frac{\cosh(n\theta/a)}{\cosh^2(n\pi b/a)}(-1)^n-a\pi \theta \su \frac{\sin((2n-1)\theta/2b)}{\sinh((2n-1)\pi a/2b)}\nonumber\\&\quad-a^2\pi^2\su \frac{\cos((2n-1)\theta/2b)\cosh((2n-1)\pi a/2b)}{\sinh^2((2n-1)\pi a/2b)}+\frac{b^2\pi^2}{2}=0,\label{3.7}\\
&b^2\pi^2\su \frac{\sinh((2n-1)\theta/2a)}{\sinh^2((2n-1)b\pi/2a)}(-1)^n+a\pi \theta \su\frac{\cos(n\theta/b)}{\cosh(n\pi a/b)}\nonumber\\
&\quad-a^2\pi^2\su\frac{\sin(n\theta/b)\sinh(n\pi a/b)}{\cosh^2(n\pi a/b)}+\frac{a\theta \pi}{2}=0.\label{3.8}
\end{align}
\end{thm}
\begin{proof} The proof is based on the functions $F_j(s,\theta;a,b)\ (j=1,2,3,4)$ and the usual residue computation. It is clear that the functions $F_i(s,\theta;a,b)\ (i=1,2,3,4)$ are meromorphic in the entire complex plane with some simple poles, see \ref{T3.1}. Here $n$ is any positive integer.

\begin{table}[htbp]\centering
\caption{Poles of functions $F_i(s,\theta;a,b)$.\label{T3.1}}
 \begin{tabular}{ll}
  \hline
  poles  & functions\\
 \hline
 $0,\ \pm n/a,\ \pm ni/b $& $F_1(s,\theta;a,b)$ \\
$\pm (2n-1)/{2a},\ \pm  (2n-1)i/{2b}$ &$F_2(s,\theta;a,b)$\\
$0,\ \pm n/a\ \pm (2n-1)i/2b $  &$F_3(s,\theta;a,b)$\\
$0,\ \pm (2n-1)/2a,\ \pm ni/b $ &$F_4(s,\theta;a,b)$\\
  \hline
 \end{tabular}
\end{table}

Since the proofs of identities \eqref{3.6}--\eqref{3.8} are similar to the proof of \eqref{3.5}, we only prove formula \eqref{3.5} here. First, we note that function $F_1(s,\theta;a,b)$ only has poles at the 0, $\pm n/a$, $\pm ni/b$ $(n\in \N)$. At an integer $\pm n/a$ the pole is simple and the residue is
\begin{align}
\underset{s=\pm n/a}\Res\Big\{ F_1(s,\theta; a,b)\Big\}= (-1)^n\frac{\pi^2}{a} \frac{\cosh(n\theta/a)}{\sinh^2(bn\pi/a)}.\label{3.9}
\end{align}
From \eqref{3.2}, at an imaginary number $\pm ni/b$, the pole has order two and the residue is
\begin{align}
\underset{s=\pm ni/b}\Res\Big\{ F_1(s,\theta; a,b)\Big\} &=\frac{\pi}{b^2}\frac{d}{ds}\left(\frac{\cosh(\theta s)}{\sin(a\pi s)}\right)|_{s=\pm ni/b}\nonumber\\
&=\frac{\theta \pi}{b^2}\frac{\sin(n\theta/b)}{\sinh(an\pi/b)}+\frac{a\pi^2}{b^2}
                \frac{\cos(n\theta/b)\cosh(an\pi/b)}{\sinh^2(an\pi/b)}.\label{3.10}
\end{align}
Furthermore, applying \eqref{e3} and \eqref{e5} we deduce the asymptotic expansion
\begin{align*}
F_1(s,\theta;a,b)=\cosh(\theta s)\left(\frac{1}{ab^2}\frac{1}{s^3}+\bigg(\frac{a}{b^2}-\frac{2}{a}\bigg)\frac{\zeta(2)}{s}+o(1)\right),\ s\rightarrow 0.
\end{align*}
Therefore, the residue of the pole of order three at 0 is found to be
\begin{align}
\underset{s=0}\Res\Big\{F_1(s,\theta; a,b)\Big\}&=\frac{1}{2}\frac{d^2}{ds^2}\left\{\cosh(\theta s)\left(\frac{1}{ab^2}+\bigg(\frac{a}{b^2}-\frac{2}{a}\bigg){\zeta(2)}s^2+o(s^3)\right)\right\}|_{s=0} \nonumber\\
&=\frac{\theta^2}{2ab^2}+\zeta(2)\Big(\frac{a}{b^2}-\frac{2}{a}\Big).\label{3.11}
\end{align}
By Lemma \ref{lem2.1}, summing the three contributions \eqref{3.9}--\eqref{3.11} yields the desired result \eqref{3.5}. Similarly, consider functions $F_j(s,\theta;a,b)\ (j=2,3,4)$. Then using residue theorem and asymptotic formulas \eqref{3.1}--\eqref{3.4}, we can deduce the evaluations \eqref{3.6}--\eqref{3.8}  by calculations similar to those above. We leave the details to the interested reader. This concludes the proof of Theorem \ref{thm3.1}.
\end{proof}

\begin{thm}\label{thm3.2}
For $a,b,\theta\in\R$ with $|\theta|<2b \pi$, $a\neq 0$ and $b>0$,
\begin{align}
&\frac{\pi^2\theta}{a^2}\su \frac{\cosh(n\theta/a)}{\sinh^2(\pi n b/a)}-2\frac{\pi^3 b}{a^2}\su \frac{\sinh(n\theta/a)\cosh(\pi n b/a)}{\sinh^3(\pi n b/a)}-\frac{\pi^2\theta}{b^2}\su \frac{\cos(n\theta/b)}{\sinh^2(\pi n a/b)}\nonumber\\
&\quad+2\frac{\pi^3 a}{b^2}\su \frac{\sin(n\theta/b)\cosh(\pi n a/b)}{\sinh^3(\pi n a/b)}+\theta\Big(\frac1{b^2}-\frac1{a^2}\Big)\zeta(2)+\frac{\theta^3}{12a^2b^2}=0,\label{3.12}\\
&\frac{\theta}{a^2}\su \frac{\cosh((2n-1)\theta/2a)}{\cosh^2((2n-1)b\pi/2a)}-2\frac{\pi b}{a^2}\su \frac{\sinh((2n-1)\theta/2a)\sinh((2n-1)b\pi/2a)}{\cosh^3((2n-1)b\pi/2a)}\nonumber\\
&\quad -\frac{\theta}{b^2}\su \frac{\cos((2n-1)\theta/2b)}{\cosh^2((2n-1)a\pi/2b)}+2\frac{\pi a}{b^2}\su \frac{\sin((2n-1)\theta/2b)\sinh((2n-1)a\pi/2b)}{\cosh^3((2n-1)a\pi/2b)}=0,\label{3.13}\\
&\frac{\theta}{a^2}\su \frac{\cosh(n\theta/a)}{\cosh^2(\pi b n/a)}-2\frac{\pi b}{a^2} \su \frac{\sinh(n\theta/a)\sinh(\pi b n/a)}{\cosh^3(\pi b n /a)}+\frac{\theta}{b^2}\su \frac{\cos((2n-1)\theta/2b)}{\sinh^2((2n-1)a\pi/2b)}\nonumber\\
&\quad -2\frac{\pi a}{b^2}\su \frac{\sin((2n-1)\theta/2b)\cosh((2n-1)a\pi/2b)}{\sinh^3((2n-1)a\pi/2b)}+\frac{\theta}{2a^2}=0.\label{3.14}
\end{align}
\end{thm}
\begin{proof}
The proofs of \eqref{3.12}--\eqref{3.14} are similar to the proof of \eqref{3.5}. We prove only \eqref{3.12} to illustrate the main ideas.

By definition of $H_1(s,\theta;a,b)$, it is easy to see that it only has poles at the $0,\ \pm n/a$ and $\pm ni/b$, where $n$ is a positive integer. Applying \eqref{3.1} and \eqref{3.2}, we find that at $\pm n/a$ and $\pm ni/b$ the poles have order two and the residues are
\begin{align}
&\underset{s=\pm n/a}\Res\Big\{H_1(s,\theta;a,b)\Big\}=\frac{\pi^2}{a^2}\left\{\frac{\theta \cosh(n\theta/a)}{\sinh^2(\pi b n/a)}-2\pi b \frac{\sinh(n\theta/a)\cosh(\pi b n/a)}{\sinh^3(\pi b n/a)}\right\},\label{3.15}\\
&\underset{s=\pm ni/b}\Res\Big\{H_1(s,\theta;a,b)\Big\}=\frac{\pi^2}{b^2}\left\{-\frac{\theta \cos(n\theta/b)}{\sinh^2(a\pi n/b)}+2a \pi \frac{\sin(n\theta/b)\cosh(a\pi n/b)}{\sinh^3(a\pi n/b)}\right\}.\label{3.16}
\end{align}
Moreover,
\begin{align*}
H_1(s,\theta;a,b)=\frac {\theta}{a^2b^2s^3}+2\zeta(2)\Big(\frac1{b^2}-\frac1{a^2}\Big)\frac{\theta }{s}+\frac{\theta^3}{6a^2b^2s}+o(1),\quad s\rightarrow 0.
\end{align*}
Hence, the residue of the pole of order three at 0 is
\begin{align}
\underset{s=0}\Res\Big\{H_1(s,\theta;a,b)\Big\}
=2\zeta(2)\theta\Big(\frac1{b^2}-\frac1{a^2}\Big)+\frac{\theta^3}{6a^2b^2}.\label{3.18}
\end{align}
Thus, by the residue theorem and using \eqref{3.15}, \eqref{3.16} and \eqref{3.18}, we obtain the desired result.

The evaluations \eqref{3.13} and \eqref{3.14} can be established similarly using $H_2(s,\theta;a,b)$ and $H_3(s,\theta;a,b)$. We leave the details to the interested reader.
\end{proof}

\subsection{Evaluations of ${\bar S}_{2p,2}(y),\ {\bar S}'_{2p-1,2}(y),\ {\bar C}_{2p,2}(y)\ {\rm and}\ {\bar C}'_{2p-1,2}(y)$}\label{sec3.2}

\begin{thm}\label{thm3.3} Let $p$ be a positive integer and $\alpha,\beta$ be real numbers such as $\alpha\beta=\pi^2$, then
\begin{align}
&\alpha^{2p} {\bar S}_{2p,2}(\alpha)-2p(-1)^{p-1}\pi^{2p-2}\beta S_{2p-1,1}(\beta)\nonumber\\&\quad\quad\quad\quad\quad\quad-(-1)^p\pi^{2p-2}\beta^2\su \frac{n^{2p}\cosh(n\beta)}{\sinh^2(n\beta)}-\delta_p=0,\label{3.19}\\
&\alpha^{2p} {\bar C}'_{2p-1,2}(\alpha)-(4p-2)(-1)^p\pi^{2p-1}C'_{2p-2,1}(\beta)\nonumber\\&\quad\quad\quad\quad\quad\quad+(-1)^p\pi^{2p-1}\beta\su \frac{(2n-1)^{2p-1}\sinh((2n-1)\beta/2)}{\cosh^2((2n-1)\beta/2)}=0,\\
&\alpha^{2p+1}{\bar C}_{2p,2}(\alpha)-(-1)^p p\frac{\pi^{2p}}{2^{2p-2}}S'_{2p-1,1}(\beta)\nonumber\\&\quad\quad\quad\quad\quad\quad+(-1)^p\frac{\pi^{2p}}{2^{2p}}\beta\su \frac{(2n-1)^{2p}\cosh((2n-1)\beta/2)}{\sinh^2((2n-1)\beta/2)}=0,\\
&\alpha^{2p}{\bar S}'_{2p-1,2}(\alpha)+2^{2p-1}(2p-1)(-1)^p\pi^{2p-1}C_{2p-2,1}(\beta)\nonumber\\&\quad\quad\quad\quad\quad\quad+2^{2p-1}(-1)^{p-1}\pi^{2p-1}\beta\su \frac{n^{2p-1}\sinh(n\beta)}{\cosh^2(n\beta)}-\chi_p=0.
\end{align}
where $\delta_1=1/2,\chi_1=\pi$ and $\delta_p=\chi_p=0$ if $p\geq 2$.
\end{thm}
\begin{proof}
Note that using the relation
\[\frac{dy}{dx}=-\frac{1}{x(1-x)z^2},\]
we can find
\begin{align}
&\su \frac{n^{2p}\cosh(ny)}{\sinh^2(ny)}=x(1-x)z^2\frac{d}{dx}S_{2p-1,1}(y),\label{3.23}\\
&\su \frac{(2n-1)^{2p-1}\sinh((2n-1)y/2)}{\cosh^2((2n-1)y/2)}=2x(1-x)z^2\frac{d}{dx}C'_{2p-2,1}(y),\label{3.24}\\
&\su \frac{(2n-1)^{2p}\cosh((2n-1)y/2)}{\sinh^2((2n-1)y/2)}=2x(1-x)\frac{d}{dx}S'_{2p-1,1}(y),\label{3.25}\\
&\su \frac{n^{2p-1}\sinh(ny)}{\cosh^2(ny)}=x(1-x)z^2\frac{d}{dx}C_{2p-2,1}(y).\label{3.26}
\end{align}
In Theorem \ref{thm3.1}, differentiating \eqref{3.5}, \eqref{3.7} $2p$ times and \eqref{3.6}, \eqref{3.8} $2p-1$ times with respect to $\theta$, respectively, then setting $\theta=0$ and letting $b\pi/a=\alpha$ and $a\pi/b=\beta$, we can produce the corresponding series with the correct powers of $n$ or $2n-1$.
Hence, from Lemma \ref{lem2.3} with $\beta=y$ and $\alpha=\pi^2/y$, Theorem \ref{thm3.3}, formulas \eqref{3.23}--\eqref{3.26} and transformation \eqref{2.23}, we may arrive at the desired identities.
\end{proof}

\begin{thm} For any $p\in\N$, the four series
\begin{align*}
{\bar S}_{2p,2}(y),\ {\bar S}'_{2p-1,2}(y),\ {\bar C}_{2p,2}(y)\ {\rm and}\ {\bar C}'_{2p-1,2}(y)
\end{align*}
can be expressed in terms of $x,z$ and $z'$.
\end{thm}

For example, setting $p=2$ and $\beta =y$ in \eqref{3.19} yields
\begin{align}\label{3.27}
{\bar S}_{4,2}(\pi^2/y)=\su \frac{n^4(-1)^{n-1}}{\sinh^2(n\pi^2/y)}=\frac{xy^5z^4}{8\pi^6}\left\{4x(1-x)yzz'+(1-x)yz^2-4\right\}.
\end{align}
Then, applying the transformation \eqref{2.23} to \eqref{3.27}, we obtain the result
\begin{align*}
{\bar S}_{4,2}(y)=\su \frac{n^4(-1)^{n-1}}{\sinh^2(ny)}=\frac{1}{8} x(1-x)z^5\big[z-4(1-x)z'\big].
\end{align*}

Next, we will provide explicit evaluations for
\[{\bar S}_{2p,2}(y),\ {\bar S}'_{2p-1,2}(y),\ {\bar C}_{2p,2}(y)\ {\rm and}\ {\bar C}'_{2p-1,2}(y)\]
with $p=1,2,3,4,5$.

\begin{cor}\label{cor3.5} Setting $r(x)=\sqrt{1-x}$, we have
\begin{align*}
&{\bar S}_{2,2}(y)= \frac{1}{8} x(1-x)z^2\big[4x(1-x)(z')^2+4(1-x)zz'-z^2\big],\\
&{\bar S}_{4,2}(y)= \frac{1}{8} x(1-x)z^5\big[z-4(1-x)z'\big],\\
&{\bar S}_{6,2}(y)= \frac{1}{8} x(1-x)z^7\big[6(1-x)(2-x)z'+(2x-3)z\big],\\
&{\bar S}_{8,2}(y)= \frac{1}{8} x(1-x)z^9\big[4(2x^3-19x^2+34x-17)z'+(3x^2-19x+17)z\big],\\
&{\bar S}_{10,2}(y)= \frac{1}{8} x(1-x)z^{11}\left\{\begin{array}{l}
10(x^4-34x^3+126x^2-155x+62)z'\\ \quad +(4x^3-102x^2+252x-155)z \end{array}\right\},\\
&{\bar S}'_{1,2}(y)= x(1-x)z^2z',\\
&{\bar S}'_{3,2}(y)=x(1-x)z^4\big[z-3(1-x)z'\big],\\
&{\bar S}'_{5,2}(y)=x(1-x)z^6\big[5(x^2-6x+5)z'+2(x-3)z\big],\\
&{\bar S}'_{7,2}(y)=x(1-x)z^8\big[7(x^3-47x^2+107x-61)z'+(3x^2-94x+107)z\big],\\
&{\bar S}'_{9,2}(y)=x(1-x)z^{10}\left\{\begin{array}{l}9(x^4-412x^3+2142x^2-3116x+1385)z'\\ \quad+4(x^3-309x^2+1071x-779)z
\end{array}\right\},\\
&{\bar C}_{2,2}(y)= \frac{r(x)}{8} xz^3\big[4(x-1)z'+z\big],\\
&{\bar C}_{4,2}(y)= \frac{r(x)}{32}xz^5\big[8(x^2-3x+2)z'+(3x-4)z\big],\\
&{\bar C}_{6,2}(y)= \frac{r(x)}{128} xz^7\big[12(x^3-17x^2+32x-16)z'+(5x^2-52x+48)z\big],\\
&{\bar C}_{8,2}(y)= \frac{r(x)}{512} xz^9\left\{\begin{array}{l}16(x^4-139x^3+546x^2-680x+272)z' \\ \quad+(7x^3-696x^2+1776x-1088) z\end{array}\right\},\\
&{\bar C}_{10,2}(y)= \frac{r(x)}{2048} xz^{11}\left\{\begin{array}{l}20(x^5-1233x^4+10400x^3-25040x^2+23808x-7936)z' \\ \quad+(9x^4-8632x^3+53232x^2-84288x+39680) z \end{array}\right\},\\
&{\bar C}'_{1,2}(y)= \frac{r(x)}{2} xz^2\big[2(x-1)z'+z\big],\\
&{\bar C}'_{3,2}(y)= -\frac{r(x)}{2} xz^4\big[6(x-1)z'+z\big],\\
&{\bar C}'_{5,2}(y)= -\frac{r(x)}{2} xz^6 \big[10(4x^2-9x+5)z'+(12x-13)z\big],\\
&{\bar C}'_{7,2}(y)= -\frac{r(x)}{2} xz^8  \big[14(16x^3-92x^2+137x-61)z'+(80x^2-292x+213)z\big],\\
&{\bar C}'_{9,2}(y)= -\frac{r(x)}{2} xz^{10}  \left\{\begin{array}{l} 18(64x^4-1168x^3+3528x^2-3809x+1385)z'\\ \quad+(448x^3-5904x^2+11688x-6233)z \end{array} \right\}.
\end{align*}
\end{cor}

\subsection{Evaluations of $S_{2p,2}(y),\ { S}'_{2p,2}(y),\ C_{2p,2}(y)\ {\rm and}\ { C}'_{2p,2}(y)$}\label{sec3.3}

In this subsection, we use Ramanujan's results of Eisenstein series to evaluate the four series
\[S_{2p,2}(y),\ { S}'_{2p,2}(y),\ C_{2p,2}(y)\ {\rm and}\ { C}'_{2p,2}(y),\]
and present explicit evaluations for $1\leq p\leq 5$. In Ramanujan's notation, the three relevant Eisenstein series are defined for $|q|<1$ by
\begin{align*}
&P(q):=1-24\su \frac{nq^n}{1-q^n},\\
&Q(q):=1+240\su \frac{n^3q^n}{1-q^n},\\
&R(q):=1-504\su \frac{n^5q^n}{1-q^n}.
\end{align*}
The functions $P,\ Q$ and $R$ were thoroughly studied in a famous paper \cite{H1927} by Ramanujan. By Entry 13 in Chapter 17 of Ramanujan's third notebook (Berndt \cite{B1991}), we have
\begin{align}
&P(q)=(1-5x)z^2+12x(1-x)zz',\label{3.52}\\
&P(q^2)=(1-2x)z^2+6x(1-x)zz',\label{3.583}\\
&Q(q)=z^4(1+14x+x^2),\label{3.54}\\
&Q(q^2)=z^4(1-x+x^2),\label{3.55}\\
&R(q)=z^6(1+x)(1-34x+x^2),\label{3.56}\\
&R(q^2)=z^6(1+x)(1-x/2)(1-2x).\label{3.57}
\end{align}
In \cite{B1991}, $L(q)=P(q),\ M(q)=Q(q)$ and $N(q)=R(q)$ in the present notation. Ramanujan \cite{H1927} proved the general result
\begin{align}\label{3.58}
\su \frac{n^pq^n}{(1-q^n)^2}=\sum C_{l,m,n}P^lQ^mR^n,
\end{align}
and derived explicit formulas for $p=2,4,6,8,10,12,14$, where $C_{l,m,n}$ is a constant and $m,n$ are non-negative integers with $l\leq 2$ and $2l+4m+6n=p+2$. Since $q=e^{-y}$, the left-hand side of Eisenstein series \eqref{3.58} can be written as
\begin{align}\label{3.59}
\su \frac{n^pq^n}{(1-q^n)^2}=\frac1{4}\su \frac{n^p}{\sinh^2(ny/2)}=\frac1{4}S_{p,2}(y/2).
\end{align}
Hence, the series $S_{2p,2}(y/2)\ (p\in\N)$ can be expressed in terms of $x,z$ and $z'$.

\begin{thm}\label{thm3.6} For all $p\in\N$, the four series
\[S_{2p,2}(y),\ { S}'_{2p,2}(y),\ C_{2p,2}(y)\ {\rm and}\ { C}'_{2p,2}(y)\]
can be expressed in terms of $x,z$ and $z'$.
\end{thm}
\begin{proof} Using the transformation \eqref{2.25} to \eqref{3.59}, we see that $S_{2p,2}(y/2)\Rightarrow S_{2p,2}(y)$, which implies that $S_{2p,2}(y)$ can be represented by $z$ and $z'$.
Then, according to Definition \ref{de1}, we deduce the relations
\begin{align}
&S'_{p,m}(y)=S_{p,m}(y/2)-2^pS_{p,m}(y),\label{3.60}\\
&{\bar C}_{p,m}(y)=C'_{p,m}(2y)-2^pC_{p,m}(2y).\label{3.61}
\end{align}
Moreover, applying the transformation \eqref{2.24} to $\sinh^2((2n-1)y/2)$, we find
\[\sinh^2((2n-1)y/2)\rightarrow -\cosh^2((2n-1)y/2).\]
Hence, the series ${ S}'_{2p,2}(y)$ and $C'_{2p,2}(y)$ can be evaluated by $z$ and $z'$.
Furthermore, using \eqref{3.61} and transformation \eqref{2.25},
we can prove that $C_{2p,2}(2y)$ can be represented by $z$ and $z'$. Finally, employing the relation
\begin{align}\label{3.62}
{C}_{p,m}(y)=C'_{p,m}(2y)+2^pC_{p,m}(2y)
\end{align}
we can now complete the proof of Theorem \ref{thm3.6}.
\end{proof}

From \cite{H1927}, we find that
\begin{align}
&\su \frac{n^2q^n}{(1-q^n)^2}=\frac{Q(q)-P^2(q)}{288},\label{3.63}\\
&\su \frac{n^4q^n}{(1-q^n)^2}=\frac{P(q)Q(q)-R(q)}{720},\label{3.64}\\
&\su \frac{n^6q^n}{(1-q^n)^2}=\frac{Q^2(q)-P(q)R(q)}{1008},\label{3.65}\\
&\su \frac{n^8q^n}{(1-q^n)^2}=\frac{P(q)Q^2(q)-Q(q)R(q)}{720},\label{3.66}\\
&\su \frac{n^{10}q^n}{(1-q^n)^2}=\frac{3Q^3(q)+2R^2(q)-5P(q)Q(q)R(q)}{1584}.\label{3.67}
\end{align}
Using \eqref{3.63}--\eqref{3.67}, we can evaluate some examples of $S_{2p,2}(y),\ { S}'_{2p,2}(y),\ C_{2p,2}(y)\ {\rm and}\ { C}'_{2p,2}(y)$ in closed form.

\begin{cor}\label{cor3.7} We have
\begin{align}
&S_{2,2}(y/2)=\frac{x(1-x)z^2}{3}\big[z^2-(1-5x)zz'-6x(1-x)(z')^2\big],\\
&S_{4,2}(y/2)=\frac{x(1-x)z^5}{30}\big[2(1+14x+x^2)z'+(7+x)z\big],\\
&S_{6,2}(y/2)=\frac{x(1-x)z^7}{42}\big[-2(1-33x-33x^2+x^3)z'+(11+22x-x^2)z\big],\\
&S_{8,2}(y/2)=\frac{x(1-x)(1+14x+x^2)z^9}{30}\big[2(1+14x+x^2)z'+(7+x)z\big],\\
&S_{10,2}(y/2)=\frac{x(1-x)z^{11}}{66}\left\{\begin{array}{l}-10(1+x)(1-20x-474x^2-20x^3+x^4)z'\\+(19+988x+1482x^2+76x^3-5x^4)z\end{array}\right\},\\
&S_{2,2}(y)=\frac{x(1-x)z^2}{24}\big[z^2-4(1-2x)zz'-12x(1-x)(z')^2\big],\\
&S_{4,2}(y)=\frac{x(1-x)z^5}{120}\big[4(1-x+x^2)z'+(2x-1)z\big],\\
&S_{6,2}(y)=\frac{x(1-x)z^7}{168}\big[2(2-x)(x+1)(2x-1)z'-(2x^2-2x-1)z\big],\\
&S_{8,2}(y)=\frac{x(1-x)(1-x+x^2)z^9}{120}\big[4(1-x+x^2)z'+(2x-1)z\big],\\
&S_{10,2}(y)=\frac{x(1-x)z^{11}}{264}\left\{\begin{array}{l} -10(2-5x+2x^2+2x^3-5x^4+2x^5)z'\\+(5-4x-6x^2+20x^3-10x^4)z \end{array}\right\},\\
&S'_{2,2}(y)=\frac{x(1-x)z^3}{6}\big[z+2(1+x)z'\big],\label{3.75}\\
&S'_{4,2}(y)=\frac{x(1-x)z^5}{30}\big[-2(7-22x+7x^2)z'+(11-7x)z\big],\label{3.76}\\
&S'_{6,2}(y)=\frac{x(1-x)z^7}{42}\big[2(31x^3-15x^2-15x+31)z'+(31x^2-10x-5)z\big],\label{3.77}\\
&S'_{8,2}(y)=\frac{x(1-x)z^9}{30}\left\{\begin{array}{l} -(254-568x+372x^2-568x^3+254x^4)z'\\ \quad+(71-93x+213x^2-127x^3)z \end{array}\right\},\label{3.78}\\
&S'_{10,2}(y)=\frac{x(1-x)z^{11}}{66}\left\{\begin{array}{l} 10(511-1261x+1006x^2+1006x^3-1261x^4+511x^5)z'\\ \quad+(-1261+2012x+3018x^2-5044x^3+2555x^4)z \end{array}\right\}.\label{b1}
\end{align}
\end{cor}
\begin{proof} These results follow from \eqref{3.52}--\eqref{3.57}, \eqref{3.59}, \eqref{3.60} and \eqref{3.63}--\eqref{3.67}.
\end{proof}

\begin{cor}\label{cor3.8}  Setting $r(x)=\sqrt{1-x}$, we have
\begin{align}
&C'_{2,2}(y)=\frac{x(1-x)z^3}{3} \big[z-(1-2x)z'\big],\label{3.79}\\
&C'_{4,2}(y)=\frac{x(1-x)z^5}{15} \big[(7+8x-8x^2)z'+2(1-2x)z\big],\label{3.80}\\
&C'_{6,2}(y)=\frac{x(1-x)z^7}{21} \big[(-31+78x-48x^2+32x^3)z'+(13-16x+16x^2)z\big],\label{3.81}\\
&C'_{8,2}(y)=\frac{x(1-x)z^9}{15} \left\{\begin{array}{l} (127-224x+96x^2+256x^3-128x^4)z'\\
\quad -4(7-6x-24x^2+16x^3)z\end{array}\right\},\label{3.82}\\
&C'_{10,2}(y)=\frac{x(1-x)z^{11}}{33} \left\{ \begin{array}{l} 5(-511+1294x-1072x^2+1568x^3-1280x^4+512x^5)z' \\ +(647-1072x+2352x^2-2560x^3+1280x^4)z\end{array} \right\},\label{b2}\\
&C'_{2,2}(2y)=\frac{xz^3}{48}\Big[\big(3r(x)+x-1\big)z-2(1-x)\big(6r(x)+x-2\big)z'\Big],\label{3.83}\\
&C'_{4,2}(2y)=\frac{xz^5}{960}\left\{\begin{array}{l} 2(1-x)\Big[-8+120r(x)+x(8-60r(x)+7x)\Big]z'\\ \quad +\Big[4-60r(x)+(3+45r(x)-7x)x\Big]z \end{array} \right\},\label{3.84}\\
&C'_{6,2}(2y)=\frac{xz^7}{5376}\left\{ \begin{array}{l} 2(1-x)\left[\begin{array}{l}\quad 32(1-63r(x))\\-x\big(48-2016r(x)+x(-78+126r(x)+31x)\big)  \end{array}\right]z' \\+ \left[\begin{array}{l} \quad 16(1-63r(x))\\-x\big(68-1092r(x)+x(-83+105r(x)+31x)\big)\end{array}\right]z  \end{array}\right\},\label{3.85}\\
&C'_{8,2}(2y)=\frac{xz^9}{15360}\left\{ \begin{array}{l} 2(1-x) \left[ \begin{array}{l} 128(-1+255r(x))+64x(4-765r(x))+127x^4\\ +48x^2(2+345r(x))-8x^3(28+15r(x))\end{array}\right]z'\\- \left[\begin{array}{l} 64(-1+255r(x))-5x^3(59+21r(x))+127x^4\\ +x\Big(16-26640r(x)+72x\big(3+145r(x)\big)\Big)\end{array}\right]z \end{array}\right\},\label{3.86}\\
&C'_{10,2}(2y)=\frac{xz^{11}}{135168}\left\{ \begin{array}{l}10(1-x)\left[ \begin{array}{l} 512(1-1023r(x))+x(-1280+1047552r(x))  \\
+x^2(1568-605088r(x))\\
+x^3(-1072+81312r(x))\\
+x^4(1294-66r(x))-511x^5\end{array}\right]z'  \\ +\left[ \begin{array}{l} (-1280+1309440r(x))+x(4416-2781504r(x)) \\ +x^2(-63252+1756656r(x))+x^3(8392-284856r(x)) \\+ x^4(-7731+297r(x))+2555x^5\end{array}\right]z \end{array} \right\}.\label{b3}
\end{align}
\end{cor}
\begin{proof}
We may apply the transformation \eqref{2.24} to \eqref{3.75}--\eqref{b1} to obtain \eqref{3.79}--\eqref{b2}, respectively. Then we may apply the transformation \eqref{2.25} to \eqref{3.79}--\eqref{b2} to obtain \eqref{3.83}--\eqref{b3}, respectively.
\end{proof}

\begin{cor}\label{cor3.9} Setting $r(x)=\sqrt{1-x}$, we have
\begin{align*}
&C_{2,2}(y)=\frac{x(1-x)z^3}{24}\Big[2(2-x)z'-z\Big],\\
&C_{4,2}(y)=\frac{x(1-x)z^5}{480}\Big[2(-8+8x+7x^2)z'+(4+7x)z\Big],\\
&C_{6,2}(y)=\frac{x(1-x)z^7}{2688}\Big[2(32-48x+78x^2-31x^3)z'-(16-52x+31x^2)z\Big],\\
&C_{8,2}(y)=\frac{x(1-x)z^9}{7680}\left\{\begin{array}{l} 2(-128+256x+96x^2-224x^3+127x^4)z'\\ \quad+(64+48x-168x^2+127x^3)z  \end{array}\right\},\\
&C_{10,2}(y)=\frac{x(1-x)z^{11}}{67584}\left\{\begin{array}{l} 10(512-1280x+1568x^2-1072x^3+1294x^4-511x^5)z'\\ -(1280-3136x+3216x^2-5176x^3+2555x^4)z \end{array} \right\},\\
&C_{2,2}(2y)=\frac{xz^3}{192} \Big[2(1-x)\big(2-x+6r(x)\big)+\big(x-1-3r(x)\big)z\Big],\\
&C_{4,2}(2y)=\frac{xz^5}{15360}\left\{\begin{array}{l} 2(1-x)\Big[-8-120r(x)+x(8+60r(x)+7x)\Big]z'\\ \quad +\Big[4+60r(x)-(-3+45r(x)+7x)x\Big]z \end{array} \right\}, \\
&C_{6,2}(2y)=\frac{xz^7}{344064}\left\{ \begin{array}{l} 2(1-x)\left[\begin{array}{l} 32(1+63r(x))-48x(1+42r(x))\\+6x^2(13+21r(x))-31x^3 \end{array}\right]z' \\+ \left[\begin{array}{l} -16(1+63r(x)) +4x(17+273r(x))\\ \quad -x^2(83+105r(x))+31x^3\end{array}\right]z  \end{array}\right\},\\
&C_{8,2}(2y)=\frac{xz^9}{3932160}\left\{ \begin{array}{l} 2(1-x) \left[ \begin{array}{l} -128(1+255r(x))+64x(4+765r(x))+127x^4\\ +48x^2(2-345r(x))+8x^3(-28+15r(x))\end{array}\right]z'\\- \left[\begin{array}{l} -64(1+255r(x))+16x(1+1665r(x))\\+72x^2(3-145r(x)) +5x^3(-59+21r(x))+127x^4 \end{array}\right]z \end{array}\right\},\\
&C_{10,2}(2y)=\frac{xz^{11}}{138412032}\left\{\begin{array}{l} 10(1-x) \left[\begin{array}{l} 512(1+1023r(x))-256x(5+4092r(x))\\+32x^2(49+18909r(x)) \\ -16x^3(67+5082r(x)) \\+2x^4(647+33r(x))-511x^5\end{array} \right]z' \\ + \left[\begin{array}{l} -1280(1+1023r(x))+192x(23+14487r(x)) \\  -16x^2(397+10979r(x))+8x^3(1049+35607r(x))\\-9x^4(859+33r(x))+2555x^5\end{array} \right]z \end{array} \right\}.
\end{align*}
\end{cor}
\begin{proof} The ten desired formulas follow without difficulty from \eqref{3.61} and \eqref{3.62} with the help of Corollaries \ref{cor3.5} and \ref{cor3.8}.
\end{proof}

\subsection{Three Differential Equations}

In Theorem \ref{thm3.2}, expanding both sides of \eqref{3.12}--\eqref{3.14} as Mclaurin series of $\theta$, then equating coefficients of $\theta^{2p+1}$ on both sides and letting $b\pi/a=\alpha$ and $a\pi/b=\beta$, and noting the elementary relations
\begin{align*}
\frac{d}{dx}\left(x^{2p+1}S_{2p,2}(x)\right)
&=(2p+1)x^{2p}\su \frac{n^{2p}}{\sinh^2(nx)}-2x^{2p+1}\su \frac{n^{2p+1}\cosh(nx)}{\sinh^3(nx)},\\
\frac{d}{dx}\left(x^{2p+1}C'_{2p,2}(x)\right)&=(2p+1)x^{2p}\su \frac{(2n-1)^{2p}}{\cosh^2((2n-1)x/2)}\\
&-x^{2p+1}\su \frac{(2n-1)^{2p+1}\sinh((2n-1)x/2)}{\cosh^3((2n-1)x/2)},
\end{align*}
and
\begin{align*}
\frac{d}{dx}\left(x^{2p+1}C_{2p,2}(x)\right)
&=(2p+1)x^{2p}\su \frac{n^{2p}}{\cosh^2(nx)}-2x^{2p+1}\su \frac{n^{2p+1}\sinh(nx)}{\cosh^3(nx)},\\
\frac{d}{dx}\left(x^{2p+1}S'_{2p,2}(x)\right)
&=(2p+1)x^{2p}\su \frac{(2n-1)^{2p}}{\sinh^2((2n-1)x/2)}\\
&-x^{2p+1}\su \frac{(2n-1)^{2p+1}\cosh((2n-1)x/2)}{\sinh^3((2n-1)x/2)},
\end{align*}
we can obtain the following three differential equations.

\begin{thm}  Let $p$ be a non-integer and $\alpha,\beta$ be real numbers such as $\alpha\beta=\pi^2$, then
\begin{align}
&\beta^{p-1}\frac{d}{d\alpha}\left(\alpha^{2p+1}S_{2p,2}(\alpha)\right)
+(-\alpha)^{p-1}\frac{d}{d\beta}\left(\beta^{2p+1}S_{2p,2}(\beta)\right)+\delta_p=0,\label{3.95}\\
&\beta^{p-1}\frac{d}{d\alpha}\left(\alpha^{2p+1}C'_{2p,2}(\alpha)\right)
+(-\alpha)^{p-1}\frac{d}{d\beta}\left(\beta^{2p+1}C'_{2p,2}(\beta)\right)=0,\label{3.96}\\
&2^{2p}\beta^{p-1}\frac{d}{d\alpha}\left(\alpha^{2p+1}C_{2p,2}(\alpha)\right)
-(-\alpha)^{p-1}\frac{d}{d\beta}\left(\beta^{2p+1}S'_{2p,2}(\beta)\right)+\lambda_p=0,\label{3.97}
\end{align}
where $\delta_0=(1/\alpha -1/\beta)/6,\ \lambda_0=1/2\beta$, $\delta_1=1/2,\ \lambda_1=0$ and $\delta_p=\lambda_p=0$ if $p\geq 2$.
\end{thm}

Setting $p=0$ in \eqref{3.95}--\eqref{3.97} and noting the fact
\[\frac{d\beta}{d\alpha}=-\frac{\beta}{\alpha},\]
we can get
\begin{align*}
&\frac{d}{d\alpha}\left(\alpha S_{0,2}(\alpha)+\beta S_{0,2}(\beta)-\frac{1}{6} (\alpha+\beta)\right)=0,\\
&\frac{d}{d\alpha}\left(\alpha C'_{0,2}(\alpha)+\beta C'_{0,2}(\beta)\right)=0,\\
&\frac{d}{d\alpha}\left(\alpha C_{0,2}(\alpha)-\beta S'_{0,2}(\beta)+\alpha/2\right)=0.
\end{align*}
Hence, by integrating we arrive at the three identities
\begin{align}
&\alpha S_{0,2}(\alpha)+\beta S_{0,2}(\beta)-\frac{1}{6} (\alpha+\beta)+1=0,\label{3.101}\\
&\alpha C'_{0,2}(\alpha)+\beta C'_{0,2}(\beta)=1,\\
&\alpha C_{0,2}(\alpha)-\beta S'_{0,2}(\beta)+\frac{\alpha}{2}=1,
\end{align}
where we have used these well-known results (\cite{B1989,B1991,L1974})
\begin{align*}
S_{0,2}(\pi)=\frac 1{6}-\frac 1{2\pi},\ C'_{0,2}(\pi)=\frac 1{2\pi}, \ C_{0,2}(\pi)=-\frac 1{2}+\frac 1{2\pi}+\frac{\Gamma^4(1/4)}{16\pi^3}
\end{align*}
and
\begin{align*}
S'_{0,2}(\pi)=-\frac 1{2\pi}+\frac{\Gamma^4(1/4)}{16\pi^3}.
\end{align*}
We point out that Berndt used the double series to prove \eqref{3.101} in his book \cite{B1989}. Please see pp. 244--245 in \cite{B1989} for the details. We also notice that $\Gamma(1/4)$ and $\pi$ (and therefore $1/\pi$) are algebraically independent (see \cite[p. 85]{BruinierGHZ2008}).

However, when $p\geq 1$, we have been unable, so far, to make any progress with the three differential equations \eqref{3.95}--\eqref{3.97}.

\section{Examples and Further Results}
In this section, we will apply the main results we have obtained in the previous sections to derive many interesting
explicit evaluations of reciprocal hyperbolic series of Ramanujan type.

\subsection{Examples}
By \eqref{1.6} if $x=1/2$, then
$$
y=\pi,\ z\Big(\frac 1{2}\Big)=\frac{\Gamma^2(1/4)}{2\pi^{3/2}},\ z'\Big(\frac 1{2}\Big)=\frac {4\sqrt{\pi}}{\Gamma^2(1/4)}.
$$
Hence, from the preceding discussion, it is easy to see that the eight series
\begin{align*}
S_{2p,2}(\pi),\ {S}'_{2p,2}(\pi),\ C_{2p,2}(\pi),\ {C}'_{2p,2}(\pi),\ {\bar S}_{2p,2}(\pi),\ {\bar S}'_{2p-1,2}(\pi),\ {\bar C}_{2p,2}(\pi)\ {\rm and}\ {\bar C}'_{2p-1,2}(\pi)
\end{align*}
can be evaluated in terms of Gamma functions and $\pi$. Moreover, by Mathematica computations, we believe that the following explicit formulas should hold.
\begin{con}\label{con-hy-ser-1} Set $\Gamma=\Gamma(1/4)$. For any $p\in\N$,
\begin{align*}
&S_{4p-2,2}(\pi)=a_p\frac{\Gamma^{8p}}{\pi^{6p}}+\eta_p,\\
&S_{4p,2}(\pi)=b_p\frac{\Gamma^{8p}}{\pi^{6p+1}},\\
&C_{4p-2,2}(\pi)=c_p\frac{\Gamma^{8p-4}}{\pi^{6p-2}}+d_p\frac{\Gamma^{8p}}{\pi^{6p}},\\
&C_{4p,2}(\pi)=e_p\frac{\Gamma^{8p}}{\pi^{6p+1}}+f_p\frac{\Gamma^{8p+4}}{\pi^{6p+3}},\\
&S'_{4p-2,2}(\pi)=g_p\frac{\Gamma^{8p-4}}{\pi^{6p-2}}+h_p\frac{\Gamma^{8p}}{\pi^{6p}},\\
&S'_{4p,2}(\pi)=i_p\frac{\Gamma^{8p}}{\pi^{6p+1}}+j_p\frac{\Gamma^{8p+4}}{\pi^{6p+3}},\\
&C'_{4p-2,2}(\pi)=k_p\frac{\Gamma^{8p}}{\pi^{6p}},\\
&C'_{4p,2}(\pi)=l_p\frac{\Gamma^{8p}}{\pi^{6p+1}},\\
&{\bar S}_{2p,2}(\pi)=m_p\frac{\Gamma^{4p}}{\pi^{3p+1}}+n_p\frac{\Gamma^{4p+4}}{\pi^{3p+3}}-\eta_p,\\
&{\bar C}_{2p,2}(\pi)=o_p\frac{\sqrt{2}\Gamma^{4p}}{\pi^{3p+1}}+p_p\frac{\sqrt{2}\Gamma^{4p+4}}{\pi^{3p+3}},\\
&{\bar S}'_{2p-1,2}(\pi)=q_p\frac{\Gamma^{4p-2}}{\pi^{3p-1/2}}+r_p\frac{\Gamma^{4p+2}}{\pi^{3p+3/2}},\\
&{\bar C}'_{2p-1,2}(\pi)=s_p\frac{\sqrt{2}\Gamma^{4p-2}}{\pi^{3p-1/2}}+t_p\frac{\sqrt{2}\Gamma^{4p+2}}{\pi^{3p+3/2}},
\end{align*}
where $\eta_1=-\frac{1}{8\pi^2},\ \eta_p=0,\ p\geq 2$ and the quantities
\[
a_p,\ b_p,\ c_p,\ d_p,\ e_p,\ f_p,\ g_p,\ h_p,\ i_p,\ j_p,\ k_p,\ l_p,\ m_p,\ n_p,\ o_p,\ p_p,\ q_p,\ r_p,\ s_p,\ t_p\in \Q.
\]
\end{con}

From Corollaries \ref{cor3.5}, \ref{cor3.7}--\ref{3.9} with the help of Mathematica, we calculated a few items of these quantities above, see Table \ref{T4.1}. The closed form of $C'_{2,2}(\pi)$ can be found in Corollary 3.9 of Berndt-Bialek-Yee's paper \cite{BB2002} or Andrews-Berndt's book \cite{AB2009}. Moreover, Table \ref{T4.1} motivates the following conjecture.

\begin{con} For any $p\in\N$, if Conjecture \ref{con-hy-ser-1} holds, then the following relations hold:
\begin{align}
c_p=\frac{g_p}{2^{4p-2}},\quad d_p=-\frac{h_p}{2^{4p-2}}=\frac{(-1)^p k_p}{2^{6p-3}},\quad e_p=-\frac{i_p}{2^{4p}}=\frac{(-1)^p l_p}{2^{6p}},\quad f_p=\frac{j_p}{2^{4p}}.
\end{align}
\end{con}

\begin{table}[htbp]\centering
\caption{Coefficients of some evaluations of reciprocal hyperbolic series.\label{T4.1}}
 \begin{tabular}{llllll}
  \hline
  $p$  &\quad 1 &\quad 2 &\quad 3 &\quad 4&\quad 5\\
 \hline\\
 $a_p $& $\frac 1{3\times 2^9}$ & $\frac 1{7\times 2^{14}}$ &$\frac{9}{11\times2^{20}}$ &$\frac{27}{2^{26}}$&$\frac{49329}{19\times2^{32}}$\\[4mm]
$b_p$ & $\frac 1{5\times 2^8}$ &$\frac{3}{5\times 2^{14}}$&$\frac{567}{65\times 2^{20}}$ & $\frac{43659}{85\times 2^{26}}$&$\frac{392931}{5\times 2^{32}}$\\[4mm]
$c_p$ & $\frac 1{2^6}$ &$\frac{9}{2^{16}}$&$\frac{945}{2^{26}}$ & $\frac{480249}{2^{36}}$&$\frac{746175969}{2^{46}}$\\[4mm]
$d_p$ &$\frac{-1}{3\times 2^9}$ &$\frac{3}{7\times 2^{19}}$&$\frac{-153}{11\times2^{29}}$ & $\frac{1701}{2^{39}}$&$\frac{-12677553}{19\times 2^{49}}$\\[4mm]
$e_p$ &$\frac{-3}{5\times 2^{11}}$ &$\frac{21}{5\times 2^{21}}$& $\frac{-18711}{65\times2^{31}}$&$\frac{5544693}{85\times 2^{41}}$ &$\frac{-201573603}{5\times 2^{51}}$\\[4mm]
$f_p$ &$\frac{1}{2^{14}}$ &$\frac{33}{2^{24}}$& $\frac{8289}{2^{34}}$ & $\frac{7660737}{2^{44}}$&$\frac{18845414721}{2^{54}}$\\[4mm]
$g_p$ &$\frac{1}{2^{4}}$ &$\frac{9}{2^{10}}$&$\frac{945}{2^{16}}$&$\frac{480249}{2^{22}}$&$\frac{746175969}{2^{28}}$\\[4mm]
$h_p$ &$\frac{1}{3\times2^{7}}$ &$\frac{-3}{7\times2^{13}}$&$\frac{153}{11\times2^{19}}$&$\frac{-1701}{2^{25}}$ &$\frac{12677553}{19\times 2^{31}}$\\[4mm]
$i_p$ &$\frac{3}{5\times2^{7}}$ &$\frac{-21}{5\times2^{13}}$&$\frac{18711}{65\times 2^{19}}$&$\frac{-5544693}{85\times2^{25}}$ &$\frac{201573603}{5\times2^{31}}$\\[4mm]
$j_p$ &$\frac{1}{2^{10}}$ &$\frac{33}{2^{16}}$&$\frac{8289}{2^{22}}$ & $\frac{7660737}{2^{28}}$&$\frac{18845414721}{2^{34}}$\\[4mm]
$k_p$ &$\frac{1}{3\times2^{6}}$ &$\frac{3}{7\times2^{10}}$& $\frac{153}{11\times 2^{14}}$&$\frac{1701}{2^{18}}$ &$\frac{12677553}{19\times2^{22}}$\\[4mm]
$l_p$ &$\frac{3}{5\times2^{5}}$ &$\frac{21}{5\times2^{9}}$&$\frac{18711}{65\times 2^{13}}$ &$\frac{5544693}{85\times2^{17}}$&$\frac{201573603}{5\times2^{21}}$\\[4mm]
$m_p$ &$\frac{1}{2^{5}}$ &$\frac{-1}{2^{7}}$&$\frac{9}{ 2^{11}}$ &$\frac{-9}{2^{11}}$&$\frac{945}{2^{17}}$\\[4mm]
$n_p$ &$\frac{-1}{2^{9}}$ &$\frac{1}{2^{11}}$&$\frac{-1}{ 2^{12}}$ &$\frac{33}{2^{17}}$&$\frac{-27}{2^{16}}$\\[4mm]
$o_p$ &$\frac{-1}{2^{5}}$ &$\frac{3}{2^{9}}$&$\frac{-99}{ 2^{15}}$ &$\frac{819}{2^{18}}$&$\frac{-171045}{2^{25}}$\\[4mm]
$p_p$ &$\frac{1}{2^{9}}$ &$\frac{-5}{2^{14}}$&$\frac{93}{ 2^{19}}$ &$\frac{-2985}{2^{24}}$&$\frac{156249}{2^{29}}$\\[4mm]
$q_p$ &$\frac{1}{2^{2}}$ &$\frac{-3}{2^{5}}$&$\frac{45}{ 2^{8}}$ &$\frac{-1071}{2^{11}}$&$\frac{44793}{2^{14}}$\\[4mm]
$r_p$ &$0$ &$\frac{1}{2^{7}}$&$\frac{-5}{ 2^{9}}$ &$\frac{243}{2^{13}}$&$\frac{-2565}{2^{14}}$\\[4mm]
$s_p$ &$\frac{-1}{2^3}$ &$\frac{3}{2^{5}}$&$\frac{-15}{ 2^{7}}$ &$\frac{189}{2^{9}}$&$\frac{-3969}{2^{11}}$\\[4mm]
$t_p$ &$\frac{1}{2^6}$ &$\frac{-1}{2^{8}}$&$\frac{7}{ 2^{10}}$ &$\frac{-87}{2^{12}}$&$\frac{1809}{2^{14}}$\\[4mm]
  \hline
 \end{tabular}
\end{table}

It is possible that some similar evaluations of infinite series involving hyperbolic functions can be established by using the methods and techniques of the present paper. For example, we can get the following examples.
\begin{exa}  Set $\Gamma=\Gamma(1/4)$. We have
\begin{align*}
&S_{2,2}\left(\frac{\pi}{2}\right)=-\frac{1}{2\pi^2}+\frac{\Gamma^4}{16\pi^4}+\frac{\Gamma^8}{192\pi^6},\\
&S_{4,2}\left(\frac{\pi}{2}\right)=\frac{11\Gamma^8}{640\pi^7}+\frac{\Gamma^{12}}{1024\pi^9},\\
&S_{6,2}\left(\frac{\pi}{2}\right)=\frac{9\Gamma^{12}}{1024\pi^{10}}+\frac{29\Gamma^{16}}{57344\pi^{12}},\\
&S_{8,2}\left(\frac{\pi}{2}\right)=\frac{363\Gamma^{16}}{40960\pi^{13}}+\frac{33\Gamma^{20}}{65536\pi^{15}},\\
&S_{10,2}\left(\frac{\pi}{2}\right)=\frac{945\Gamma^{20}}{65536\pi^{16}}+\frac{4761\Gamma^{24}}{5767168\pi^{18}},\\
&S_{12,2}\left(\frac{\pi}{2}\right)=\frac{1179927\Gamma^{24}}{34078720\pi^{19}}+\frac{8289\Gamma^{24}}{4194304\pi^{21}},\\
&C_{2,2}(2\pi)=\frac{(2\sqrt{2}+1)\Gamma^4}{512\pi^4}-\frac{(3\sqrt{2}+1)\Gamma^8}{12288\pi^6},\\
&C_{4,2}(2\pi)=-\frac{3(20\sqrt{2}+1)\Gamma^8}{327680\pi^7}+\frac{(5\sqrt{2}+1)\Gamma^{12}}{524288\pi^9},\\
&C_{6,2}(2\pi)=\frac{9(22\sqrt{2}+1)\Gamma^{12}}{8388608\pi^{10}}-\frac{3(217\sqrt{2}-1)\Gamma^{16}}{469762048\pi^{12}},\\
&C_{8,2}(2\pi)=-\frac{21(1560\sqrt{2}-1)\Gamma^{16}}{5368709120\pi^{13}}+\frac{3(995\sqrt{2}+11)\Gamma^{20}}{8589934592\pi^{15}},\\
&C_{10,2}(2\pi)=\frac{945(362\sqrt{2}+1)\Gamma^{20}}{137438953472\pi^{16}}+\frac{9(190971\sqrt{2}+17)\Gamma^{24}}{12094627905536\pi^{18}},\\
&C'_{2,2}(2\pi)=-\frac{(2\sqrt{2}-1)\Gamma^4}{128\pi^4}+\frac{(3\sqrt{2}-1)\Gamma^8}{3072\pi^6},\\
&C'_{4,2}(2\pi)=\frac{3(20\sqrt{2}-1)\Gamma^8}{20480\pi^7}-\frac{(5\sqrt{2}-1)\Gamma^{12}}{32768\pi^9},\\
&C'_{6,2}(2\pi)=-\frac{9(22\sqrt{2}-1)\Gamma^{12}}{131072\pi^{10}}+\frac{3(217\sqrt{2}+1)\Gamma^{16}}{7340032\pi^{12}},\\
&C'_{8,2}(2\pi)=\frac{21(1560\sqrt{2}+1)\Gamma^{16}}{20971520\pi^{13}}-\frac{3(995\sqrt{2}-11)\Gamma^{20}}{33554432\pi^{15}}\\
&C'_{10,2}(2\pi)=-\frac{945(362\sqrt{2}-1)\Gamma^{20}}{134217728\pi^{16}}+\frac{9(190971\sqrt{2}-17)\Gamma^{24}}{1181160064\pi^{18}}.
\end{align*}
\end{exa}

\begin{exa}  Set $\Gamma=\Gamma(1/4)$. We have
\begin{align*}
&\su \frac{n^2\cosh(n\pi)}{\sinh^2(n\pi)}=-\frac{1}{8\pi^2}+\frac{\Gamma^4}{32\pi^4}+\frac{\Gamma^8}{512\pi^6},\\
&\su \frac{n^4\cosh(n\pi)}{\sinh^2(n\pi)}=\frac{\Gamma^8}{128\pi^7}+\frac{\Gamma^{12}}{2048\pi^9},\\
&\su \frac{n^6\cosh(n\pi)}{\sinh^2(n\pi)}=\frac{9\Gamma^{12}}{2^{11}\pi^{10}}+\frac{\Gamma^{16}}{2^{12}\pi^{12}},\\
&\su \frac{n^8\cosh(n\pi)}{\sinh^2(n\pi)}=\frac{9\Gamma^{16}}{2^{11}\pi^{13}}+\frac{33\Gamma^{16}}{2^{17}\pi^{15}},\\
&\su \frac{n^{10}\cosh(n\pi)}{\sinh^2(n\pi)}=\frac{945\Gamma^{20}}{2^{17}\pi^{16}}+\frac{27\Gamma^{24}}{2^{16}\pi^{18}},\\
&\su \frac{(2n-1)^2\cosh((2n-1)/2)}{\sinh^2((2n-1)\pi/2)}=\frac{\sqrt{2}\Gamma^4}{8\pi^4}+\frac{\sqrt{2}\Gamma^8}{128\pi^6},\\
&\su \frac{(2n-1)^4\cosh((2n-1)/2)}{\sinh^2((2n-1)\pi/2)}=\frac{3\sqrt{2}\Gamma^8}{32\pi^7}+\frac{5\sqrt{2}\Gamma^{12}}{1024\pi^9},\\
&\su \frac{(2n-1)^6\cosh((2n-1)/2)}{\sinh^2((2n-1)\pi/2)}=\frac{99\sqrt{2}\Gamma^{12}}{2^9\pi^{10}}+\frac{93\sqrt{2}\Gamma^{16}}{2^{13}\pi^{12}},\\
&\su \frac{(2n-1)^8\cosh((2n-1)/2)}{\sinh^2((2n-1)\pi/2)}=\frac{819\sqrt{2}\Gamma^{16}}{2^{10}\pi^{13}}+\frac{2985\sqrt{2}\Gamma^{20}}{2^{16}\pi^{15}},\\
&\su \frac{(2n-1)^{10}\cosh((2n-1)/2)}{\sinh^2((2n-1)\pi/2)}=\frac{171045\sqrt{2}\Gamma^{20}}{2^{15}\pi^{16}}+\frac{156249\sqrt{2}\Gamma^{24}}{2^{19}\pi^{18}}.
\end{align*}
\end{exa}

\begin{exa}  Set $\Gamma=\Gamma(1/4)$. We have
\begin{align*}
&\su \frac{n^3\cosh(n\pi)}{\sinh^3(n\pi)}(-1)^{n-1}=\frac 1{16\pi^3}+\frac{3\Gamma^4}{128\pi^5}-\frac{3\Gamma^8}{1024\pi^7}+\frac{\Gamma^{12}}{8192\pi^9},\\
&\su \frac{n^5\cosh(n\pi)}{\sinh^3(n\pi)}(-1)^{n-1}=-\frac{5\Gamma^8}{512\pi^8}+\frac{5\Gamma^{12}}{4096\pi^{10}}-\frac{\Gamma^{16}}{32768\pi^{12}},\\
&\su \frac{n^7\cosh(n\pi)}{\sinh^3(n\pi)}(-1)^{n-1}=\frac{63\Gamma^{12}}{8192\pi^{11}}-\frac{7\Gamma^{16}}{8192\pi^{13}}+\frac{13\Gamma^{20}}{524288\pi^{15}},\\
&\su \frac{n^9\cosh(n\pi)}{\sinh^3(n\pi)}(-1)^{n-1}=-\frac{81\Gamma^{16}}{8192\pi^{14}}+\frac{297\Gamma^{20}}{262144\pi^{16}}-\frac{17\Gamma^{24}}{524288\pi^{18}},\\
&\su \frac{n^{11}\cosh(n\pi)}{\sinh^3(n\pi)}(-1)^{n-1}=\frac{10395\Gamma^{20}}{524288\pi^{17}}-\frac{297\Gamma^{24}}{131072\pi^{19}}+\frac{2169\Gamma^{28}}{33554432\pi^{21}},\\
&\su \frac{n^3\sinh(n\pi)}{\cosh^3(n\pi)}(-1)^{n-1}=-\frac{3\sqrt{2}\Gamma^4}{2^7\pi^5}+\frac{3\sqrt{2}\Gamma^8}{2^{10}\pi^7}-\frac{\sqrt{2}\Gamma^{12}}{2^{14}\pi^9},\\
&\su \frac{n^5\sinh(n\pi)}{\cosh^3(n\pi)}(-1)^{n-1}=\frac{15\sqrt{2}\Gamma^8}{2^{11}\pi^8}-\frac{25\sqrt{2}\Gamma^{12}}{2^{15}\pi^{10}}+\frac{13\sqrt{2}\Gamma^{16}}{2^{19}\pi^{12}},\\
&\su \frac{n^7\sinh(n\pi)}{\cosh^3(n\pi)}(-1)^{n-1}=-\frac{693\sqrt{2}\Gamma^{12}}{2^{17}\pi^{11}}+\frac{65\sqrt{2}\Gamma^{16}}{2^{20}\pi^{13}}-\frac{293\sqrt{2}\Gamma^{20}}
{2^{24}\pi^{15}},\\
&\su \frac{(2n-1)^2\cosh((2n-1)\pi/2)}{\sinh^3((2n-1)\pi/2)}(-1)^{n-1}=\frac{\Gamma^2}{4\pi^{7/2}}+\frac{\Gamma^{10}}{512\pi^{15/2}},\\
&\su \frac{(2n-1)^2\sinh((2n-1)\pi/2)}{\cosh^3((2n-1)\pi/2)}(-1)^{n-1}=-\frac{\sqrt{2}\Gamma^2}{2^3\pi^{7/2}}+\frac{\sqrt{2}\Gamma^6}{2^5\pi^{11/2}},\\
&\su \frac{(2n-1)^4\sinh((2n-1)\pi/2)}{\cosh^3((2n-1)\pi/2)}(-1)^{n-1}=\frac{3\sqrt{2}\Gamma^6}{2^4\pi^{13/2}}-\frac{\sqrt{2}\Gamma^{10}}{2^6\pi^{17/2}}
+\frac{\sqrt{2}\Gamma^{14}}{2^{11}\pi^{21/2}},\\
&\su \frac{(2n-1)^6\sinh((2n-1)\pi/2)}{\cosh^3((2n-1)\pi/2)}(-1)^{n-1}=-\frac{45\sqrt{2}\Gamma^{10}}{2^7\pi^{19/2}}+\frac{21\sqrt{2}\Gamma^{14}}{2^9\pi^{23/2}}
-\frac{5\sqrt{2}\Gamma^{18}}{2^{12}\pi^{27/2}},\\
&\su \frac{(2n-1)^8\sinh((2n-1)\pi/2)}{\cosh^3((2n-1)\pi/2)}(-1)^{n-1}=\frac{189\sqrt{2}\Gamma^{14}}{2^7\pi^{25/2}}-\frac{87\sqrt{2}\Gamma^{18}}{2^9\pi^{29/2}}
+\frac{157\sqrt{2}\Gamma^{22}}{2^{15}\pi^{33/2}},\\
&\su \frac{(2n-1)^{10}\sinh((2n-1)\pi/2)}{\cosh^3((2n-1)\pi/2)}(-1)^{n-1}=-\frac{19845\sqrt{2}\Gamma^{18}}{2^{11}\pi^{31/2}}+\frac{9045\sqrt{2}\Gamma^{22}}{2^{13}\pi^{35/2}}
-\frac{1035\sqrt{2}\Gamma^{26}}{2^{15}\pi^{39/2}}.
\end{align*}
\end{exa}

\subsection{Further Results}

From \cite{L1974}, we find
\begin{align}
&\frac1{(2s-1)!}\frac{d^{2s-2}}{dy^{2s-2}}\left(\frac{1}{\sinh^2(y)}\right)=\sum\limits_{k=1}^s \frac{A_{2s,2k}}{\sinh^{2k}(y)},\label{4.14}\\
&\frac1{(2s-1)!}\frac{d^{2s-2}}{dy^{2s-2}}\left(\frac{1}{\cosh^2(y)}\right)=\sum\limits_{k=1}^s (-1)^{k+1} \frac{A_{2s,2k}}{\cosh^{2k}(y)},\label{4.15}
\end{align}
where the coefficients $A_{2s,2k} \ (A_{2s,2s}:=1)$ have the following recurrence relation:
\begin{align*}
A_{2s+2,2k}=\frac{1}{2s(2s+2)}\left((2k-1)(2k-2)A_{2s,2k-2}+4k^2A_{2s,2k}\right)
\end{align*}
for $1\leq k\leq s$. See Table 5 of \cite{L1974} for some values of $A_{2s,2k}$. For $s\geq 1$ and $p\geq 0$,
using notations \eqref{4.14} and \eqref{4.15} one obtains
\begin{align}
&\frac{1}{(2s-1)!}\frac{d^{2s-2}}{dy^{2s-2}}S_{p,2}(y)=\sum\limits_{k=1}^s A_{2s,2k}S_{p+2s-2,2k}(y),\label{4.16}\\
&\frac{1}{(2s-1)!}\frac{d^{2s-2}}{dy^{2s-2}}{\bar S}_{p,2}(y)=\sum\limits_{k=1}^s A_{2s,2k}{\bar S}_{p+2s-2,2k}(y),\label{4.17}\\
&\frac{1}{(2s-1)!}\frac{d^{2s-2}}{dy^{2s-2}}C_{p,2}(y)=\sum\limits_{k=1}^s (-1)^{k+1}A_{2s,2k}C_{p+2s-2,2k}(y),\label{4.18}\\
&\frac{1}{(2s-1)!}\frac{d^{2s-2}}{dy^{2s-2}}{\bar C}_{p,2}(y)=\sum\limits_{k=1}^s (-1)^{k+1}A_{2s,2k}{\bar C}_{p+2s-2,2k}(y),\label{4.19}\\
&\frac{1}{(2s-1)!}\frac{d^{2s-2}}{dy^{2s-2}}S'_{p,2}(y)=\frac{1}{2^{2s-2}}\sum\limits_{k=1}^s A_{2s,2k}S'_{p+2s-2,2k}(y),\label{4.20}\\
&\frac{1}{(2s-1)!}\frac{d^{2s-2}}{dy^{2s-2}}{\bar S'}_{p,2}(y)=\frac{1}{2^{2s-2}}\sum\limits_{k=1}^s A_{2s,2k}{\bar S'}_{p+2s-2,2k}(y),\label{4.21}\\
&\frac{1}{(2s-1)!}\frac{d^{2s-2}}{dy^{2s-2}}C'_{p,2}(y)=\frac{1}{2^{2s-2}}\sum\limits_{k=1}^s (-1)^{k+1}A_{2s,2k}C'_{p+2s-2,2k}(y),\label{4.22}\\
&\frac{1}{(2s-1)!}\frac{d^{2s-2}}{dy^{2s-2}}{\bar C'}_{p,2}(y)=\frac{1}{2^{2s-2}}\sum\limits_{k=1}^s (-1)^{k+1} A_{2s,2k}{\bar C'}_{p+2s-2,2k}(y).\label{4.23}
\end{align}
From the recurrence relations \eqref{4.16}--\eqref{4.23}, we find
\begin{align*}
&S_{p+2k-2,2k}(y)\in \Q\left[ \frac{d^{2k-2}}{dy^{2k-2}}S_{p,2}(y),\frac{d^{2k-4}}{dy^{2k-4}}S_{p+2,2}(y),\cdots,\frac{d^{2}}{dy^{2}}S_{p+2k-4,2}(y),S_{p+2k-2,2}(y)\right],\\
&{\bar S}_{p+2k-2,2k}(y)\in \Q\left[ \frac{d^{2k-2}}{dy^{2k-2}}{\bar S}_{p,2}(y),\frac{d^{2k-4}}{dy^{2k-4}}{\bar S}_{p+2,2}(y),\cdots,\frac{d^{2}}{dy^{2}}{\bar S}_{p+2k-4,2}(y),{\bar S}_{p+2k-2,2}(y)\right],\\
&C_{p+2k-2,2k}(y)\in \Q\left[ \frac{d^{2k-2}}{dy^{2k-2}}C_{p,2}(y),\frac{d^{2k-4}}{dy^{2k-4}}C_{p+2,2}(y),\cdots,\frac{d^{2}}{dy^{2}}C_{p+2k-4,2}(y),C_{p+2k-2,2}(y)\right],\\
&{\bar C}_{p+2k-2,2k}(y)\in \Q\left[ \frac{d^{2k-2}}{dy^{2k-2}}{\bar C}_{p,2}(y),\frac{d^{2k-4}}{dy^{2k-4}}{\bar C}_{p+2,2}(y),\cdots,\frac{d^{2}}{dy^{2}}{\bar C}_{p+2k-4,2}(y),{\bar C}_{p+2k-2,2}(y)\right],\\
&S'_{p+2k-2,2k}(y)\in \Q\left[ \frac{d^{2k-2}}{dy^{2k-2}}S'_{p,2}(y),\frac{d^{2k-4}}{dy^{2k-4}}S'_{p+2,2}(y),\cdots,\frac{d^{2}}{dy^{2}}S'_{p+2k-4,2}(y),S'_{p+2k-2,2}(y)\right],\\
& {\bar S'}_{p+2k-2,2k}(y)\in \Q\left[ \frac{d^{2k-2}}{dy^{2k-2}}{\bar S'}_{p,2}(y),\frac{d^{2k-4}}{dy^{2k-4}}{\bar S'}_{p+2,2}(y),\cdots,\frac{d^{2}}{dy^{2}}{\bar S'}_{p+2k-4,2}(y),{\bar S'}_{p+2k-2,2}(y)\right],\\
& C'_{p+2k-2,2k}(y)\in \Q\left[ \frac{d^{2k-2}}{dy^{2k-2}}C'_{p,2}(y),\frac{d^{2k-4}}{dy^{2k-4}}C'_{p+2,2}(y),\cdots,\frac{d^{2}}{dy^{2}}C'_{p+2k-4,2}(y),C'_{p+2k-2,2}(y)\right],\\
& {\bar C'}_{p+2k-2,2k}(y)\in \Q\left[ \frac{d^{2k-2}}{dy^{2k-2}}{\bar C'}_{p,2}(y),\frac{d^{2k-4}}{dy^{2k-4}}{\bar C'}_{p+2,2}(y),\cdots,\frac{d^{2}}{dy^{2}}{\bar C'}_{p+2k-4,2}(y),{\bar C'}_{p+2k-2,2}(y)\right].
\end{align*}
Hence, all series of the forms ($p\geq k$)
\begin{align*}
S_{2p,2k}(y),\ {\bar S}_{2p,2k}(y),\ C_{2p,2k}(y),\ {\bar C}_{2p,2k}(y),\ S'_{2p,2k}(y),\ {\bar S'}_{2p-1,2k}(y),\ C'_{2p,2k}(y),\ {\bar C'}_{2p-1,2k}(y)
\end{align*}
can be evaluated. So, when $x=1/2,y=\pi$, these series can be evaluated by special values of the Gamma function and $\pi$. Next, we evaluate some specific cases with the help of Mathematica.  Set $\Gamma=\Gamma(1/4)$. Then

\begin{align*}
S_{4,4}(\pi)=\su \frac{n^4}{\sinh^4(n\pi)}=&-\frac{1}{32\pi^4}-\frac{\Gamma^8}{1920\pi^7}+\frac{\Gamma^8}{1024\pi^8}+\frac{\Gamma^{16}}{393216\pi^{12}},\\
{\bar S}_{4,4}(\pi)=\su \frac{n^4(-1)^{n-1}}{\sinh^4(n\pi)}=&\frac{1}{32\pi^4}+\frac{\Gamma^4}{64\pi^6}+\frac{\Gamma^8}{192\pi^7}-\frac{3\Gamma^8}{1024\pi^8}-\frac{\Gamma^{12}}{3072\pi^9}
 +\frac{\Gamma^{12}}{4096\pi^{10}}-\frac{\Gamma^{16}}{393216\pi^{12}},\\
C_{4,4}(\pi)=\su \frac{n^4}{\cosh^4(n\pi)}=&-\frac{\Gamma^4}{128\pi^6}-\frac{\Gamma^8}{5120\pi^7}+\frac{\Gamma^8}{1024\pi^8}+\frac{\Gamma^{12}}{24576\pi^9}
-\frac{\Gamma^{12}}{8192\pi^{10}},\\
{\bar C}_{4,4}(\pi)=\su \frac{n^4(-1)^{n-1}}{\cosh^4(n\pi)}=&\frac{\sqrt{2}\Gamma^4}{64\pi^6}+\frac{\sqrt{2}\Gamma^8}{256\pi^7}-\frac{3\sqrt{2}\Gamma^8}{1024\pi^8}
-\frac{5\sqrt{2}\Gamma^{12}}{24576\pi^9} +\frac{\sqrt{2}\Gamma^{12}}{8192\pi^{10}}-\frac{\sqrt{2}\Gamma^{16}}{262144\pi^{12}}
\end{align*}
and
\begin{align*}
S'_{4,4}(\pi)=\su \frac{(2n-1)^4}{\sinh^4((2n-1)\pi/2)}=&\frac{\Gamma^4}{8\pi^6}-\frac{\Gamma^8}{320\pi^7}+\frac{\Gamma^8}{64\pi^8}-\frac{\Gamma^{12}}{1536\pi^9}
+\frac{\Gamma^{12}}{512\pi^{10}},\\
C'_{4,4}(\pi)=\su \frac{(2n-1)^4}{\cosh^4((2n-1)\pi/2)}=& \frac{\Gamma^8}{80\pi^7}-\frac{\Gamma^8}{32\pi^8},\\
{\bar S'}_{3,4}(\pi)=\su \frac{(2n-1)^3(-1)^{n-1}}{\sinh^4((2n-1)\pi/2)}=&\frac{\Gamma^2}{4\pi^{9/2}}+\frac{\Gamma^6}{16\pi^{11/2}}-\frac{\Gamma^{10}}{192\pi^{15/2}}
+\frac{3\Gamma^{10}}{512\pi^{17/2}},\\
{\bar C'}_{3,4}(\pi)=\su \frac{(2n-1)^3(-1)^{n-1}}{\cosh^4((2n-1)\pi/2)}=&\frac{\sqrt{2}\Gamma^2}{8\pi^{9/2}}+\frac{\sqrt{2}\Gamma^6}{16\pi^{11/2}}
-\frac{3\sqrt{2}\Gamma^6}{64\pi^{13/2}}-\frac{\sqrt{2}\Gamma^{10}}{384\pi^{15/2}}.
\end{align*}

These results can be checked numerically by Mathematica or Maple. It is highly likely that the closed form representations of some other infinite series involving hyperbolic functions can be proved by using the techniques of this paper.

\section{Concluding Remarks}
In this paper we have evaluated eight families of infinite series involving hyperbolic functions which include some identities discovered by Ramanujan as special cases. The key idea is based on their contour integral representations and residue computations with the help of some well-known results of Eisenstein series given by Ramanujan, Berndt et al. These explicit evaluations have played important roles in the computation of several classes infinite integrals of Berndt type (see \cite{RXYZ2024,XuZhao2024Jan}).

Our method works well when the parameter $p>0$ in all the eight families of the reciprocal hyperbolic series but does not apply in general to $p\le 0$ cases even though the series still converge. Therefore, new ideas and techniques are required to fully reveal the structure of these series in general.

\medskip
\noindent
\textbf{Acknowledgement}. Ce Xu is supported by the National Natural Science Foundation of China (Grant No. 12101008), the Natural Science Foundation of Anhui Province (Grant No. 2108085QA01) and the University Natural Science Research Project of Anhui Province (Grant No. KJ2020A0057). Jianqiang Zhao is supported by the Jacobs Prize from The Bishop's School.


\begin{thebibliography}{999}

\bibitem{A2000}
Andrews, G.E.; Askey, R.; Roy, R. \emph{Special functions};
Cambridge University Press, 2000; pp. 146--149.

\bibitem{AB2009}
Andrews, G.E.; Berndt, B.C. \emph{Ramanujan's lost notebook, Part II}; Springer, New York, 2009; pp. 223--233.

\bibitem{AB2013}
Andrews, G.E.; Berndt, B.C. \emph{Ramanujan's lost notebook, Part IV}; Springer, New York, 2013; pp. 265--284.

\bibitem{AskeyOD1999}
Askey, R.A.; Olde Daalhuis, A.B. Generalized Hypergeometric Functions and Meijer
$G$-Function, in: \emph{NIST Handbook of Mathematical Functions}; Olver, Lozier, Boisvert and Clark
Eds.;  NIST and Cambridge University Press; Cambridge, New York, Melbourne, Madrid, Cape Town, Singapore, S\~{a}o Paulo, Delhi, Dubai, Tokyo, 2010; pp. 403--418.

\bibitem{B1977}
Berndt, B.C. Modular transformations ang generalizations of several formulae of Ramanujan, \emph{Rocky Mt. J. Math.} \textbf{1977}, \emph{7}, 147--189.

\bibitem{B1978}
Berndt, B.C. Analytic Eisenstein series, theta-functions, and series relations in the spirit of Ramanujan,
\emph{J. Reine Angew. Math.} \textbf{1978} \textbf{304}, 332--365.

\bibitem{B1989}
Berndt, B.C. \emph{Ramanujan's notebooks, part II}; Springer-Verlag, New York, 1989; pp. 241--299.

\bibitem{B1991}
Berndt, B.C. \emph{Ramanujan's notebooks, part III};, Springer-Verlag, New York, 1990; pp. 87--142.

\bibitem{B2004}
Berndt, B.C. An unpublished manuscript of Ramanujan on infinite series identities, \emph{J. Ramanujan
Math. Soc.} \textbf{2004}, \emph{19}, 57--74.

\bibitem{B2016}
Berndt, B.C. Integrals associated with Ramanujan and elliptic functions, \emph{Ramanujan J.} \textbf{2016}, \emph{41}, 369--389.

\bibitem{BB2002}
Berndt, B.C.; Bialek, P.R.;  Yee, A.J.\ Formulas of Ramanujan for the power series coefficients of certain quotients of Eisenstein series, \emph{Int. Math. Res. Notices} \textbf{2002}, \emph{21}, 1077--1109.

\bibitem{BS2017}
Berndt, B.C.; Straub, A.\ \emph{Ramanujan's formula for $\zeta(2n+1)$}, Springer, Cham, 2017; pp. 13--34.

\bibitem{BruinierGHZ2008}
Bruinier, J.H.; van der Geer, G.; Harder, G.; Zagier, D.\ \emph{The 1-2-3 of Modular Forms}, Lectures at a Summer School in Nordfjordeid, Norway (Universitext), Springer-Verlag, Berlin, Heidelberg, 2008.

\bibitem{C1889}
Cauchy, A.\ \emph{Oeuvres, S$\acute{\rm e}$rie II, t. VII}; Gauthier-Villars, Paris, 1889.

\bibitem{C1961}
Cayley, A.\ \emph{An elementary treatise on Elliptic functions}, 2nd ed; Dover, New York, 1961; pp. 245--301.

\bibitem{FS1998}
Flajolet, P.;  Salvy, B.\ Euler sums and contour integral representations, \emph{Exp. Math.} \textbf{1998}, \emph{7}, 15--35.

\bibitem{G1892}
Greenhill, A.G.\ \emph{The applications of Eilliptic functions}; Macmillan, London, 1892; pp. 284--286.

\bibitem{G1989}
Guo, D.R.; Tu, F.-Y.; Wang, Z.X.\ \emph{Special functions}; World Scientific Publishing Co Pte Ltd, 1989; pp. 498--574.

\bibitem{H1958}
Hancock, H.\ \emph{Theory of elliptic functions}; Dover, New York (1958).

\bibitem{H1927}
Hardy, G.H.; Seshu Aiyar, P.V.; Wilson, B.M.\ \emph{Collected papers of Srinivasa Ramanujan}; Cambridge University Press, 1927; pp. 134--147.

\bibitem{T2015}
Komori, Y.; Matsumoto, K.; Tsumura, H.\ Infinite series involving hyperbolic functions, \emph{Lith. Math. J.} \textbf{2015}, \emph{55}, 102--118.

\bibitem{L2011}
Lim, S.G.\ Identities about infinite series containing hyperbolic functions and trigonometric functions, \emph{Korean J. Math.} \textbf{2011}, \emph{19}, 465--480.

\bibitem{L1974}
Ling, C.B. On summation of series of hyperbolic functions, \emph{Siam J. Math. Anal.} \textbf{1974}, \emph{5}, 551--562.

\bibitem{L1975}
Ling, C.B. On summation of series of hyperbolic functions, II, \emph{Siam J. Math. Anal.} \textbf{1975}, \emph{6}, 129--139.

\bibitem{L1978}
Ling, C.B. Generalization  of certain  summations  due  to Ramanujan, \emph{Siam J. Math. Anal.} \textbf{1978}, \emph{9}, 34--48.

\bibitem{L1988}
Ling, C.B. Extension of summation of series of hyperbolic functions, \emph{J. Math. Anal. Appl.} \textbf{1988}, \emph{131}, 451--464.

\bibitem{R2012}
Ramanujan, S.\ \emph{Notebooks (2 volumes)}. Tata institute of fundamental research; Bombay 1957; second ed., 2012.

\bibitem{R1988}
Ramanujan, S.\ \emph{The lost notebook and other unpublished papers}; Narosa, New Delhi; 1988.

\bibitem{RXYZ2024}
Rui, H.; Xu, C.; Yang, R.; Zhao, J. Explicit evaluation of a family of Berndt-type integrals,
\emph{J.\ Math.\ Analysis Appl.} \textbf{2025}, \emph{541} (1), Art. 128676. doi.org/10.1016/j.jmaa.2024.128676.

\bibitem{S1954}
Sandham, H.F.\ Some infinite series, \emph{Proc.\ London Math.\ Soc.} \textbf{1954}, \emph{5}, 430--436.

\bibitem{T2008}
Tsumura, H.\ On certain analogues of Eisenstein series and their evaluation formulas of Hurwitz type, \emph{Bull.\ London Math.\ Soc.} \textbf{2008}, \emph{40}, 289--297.

\bibitem{T2010}
Tsumura, H.\ Analogues of the Hurwitz formulas for level 2 Eisenstein series, \emph{Results in Math.} \textbf{2010}, \emph{58}, 365--378.

\bibitem{T2012}
Tsumura, H.\ Analogues of level-$N$ Eisenstein series, \emph{Pacific J.\ Math.} \textbf{2012}, \emph{265}, 489--510.

\bibitem{WW1966}
Whittaker, E.T.;  Watson, G.N.\ \emph{A course of modern analysis}, 4th ed.; Cambridge University Press, Cambridge, 1966; pp. 511--512.

\bibitem{X2018}
Xu, C.\ Some evaluation of infinite series involving trigonometric and hyperbolic Functions, \emph{Results in Math.} \textbf{2018}, \emph{73}, 1--18.

\bibitem{XuZhao2024Jan}
Xu, C.; Zhao, J. General Berndt-type integrals and series associated with Jacobi elliptic functions. arXiv:2401.01385.

\bibitem{Ya-2018}
Yakubovich, S.\ On the curious series related to the elliptic integrals, \emph{Ramanujan J.} \textbf{2018}, \emph{45}, 797--815.

\bibitem{Z1979}
Zucker, I.J.\ The summation of series of hyperbolic functions, \emph{Siam J.\ Math.\ Anal.} \textbf{1979}, \emph{10}, 192--206.

\bibitem{Z1984}
Zucker, I.J.\ Some infinite series of exponential and hyperbolic functions, \emph{Siam J. Math. Anal.} \textbf{1984}, \emph{15}, 406--413.

\end{thebibliography}
\end{document}